\documentclass[fontsize=12pt,a4paper,headings=normal,
twoside=false,leqno,parskip=half-,abstract=true]{scrartcl}
\usepackage[english]{babel}
\usepackage[utf8]{inputenc}
\usepackage{hyperref}
\hypersetup{
 pdftitle={Conley},
 pdfauthor={Phillipo Lappicy},
 colorlinks=true,
 linkcolor=blue,
 citecolor=blue,
 filecolor=blue,
 urlcolor=blue}

\usepackage{graphicx}
\usepackage[format=plain,labelfont=bf,font=small]{caption}
\usepackage{xcolor}
\usepackage[arrow, matrix, curve]{xy}
\usepackage{float}

\usepackage{caption}
\usepackage{tabulary}
\usepackage{array}
\newcolumntype{N}[1]{>{\centering\arraybackslash}m{#1}}

\usepackage{amsmath,amsthm}
\usepackage{amssymb} 
\usepackage{blkarray,bigstrut}

\usepackage{tikz-cd}

\usepackage{pgfplots}
\pgfplotsset{compat=1.13}

\usetikzlibrary{matrix,arrows,positioning,calc}

\tikzset{
  curvedlink/.style={
    to path={
      let \p1=(\tikztostart.east), \p2=(\tikztotarget.west),
      \n1= {abs(\y2-\y1)/4} in
      (\p1) arc(90:-90:\n1) -- ([yshift=2*\n1]\p2) arc (90:270:\n1)
    },
  }
}

\makeatletter
\newcommand{\tpitchfork}{%
  \vbox{
    \baselineskip\z@skip
    \lineskip-.52ex
    \lineskiplimit\maxdimen
    \m@th
    \ialign{##\crcr\hidewidth\smash{$-$}\hidewidth\crcr$\pitchfork$\crcr}
  }%
}
\makeatother
\usepackage{latexsym}

\usepackage[notref,notcite,color,final 
]{showkeys}

\definecolor{refkey}{rgb}{1,0,0}
\definecolor{labelkey}{rgb}{1,0,0}

\usepackage{tikz}

\usepackage[textwidth=2cm,textsize=small,backgroundcolor=none]{todonotes}

  \mathchardef\ordinarycolon\mathcode`\:
  \mathcode`\:=\string"8000
  \begingroup \catcode`\:=\active
    \gdef:{\mathrel{\mathop\ordinarycolon}}
  \endgroup

\newtheorem{thm}{Theorem}[section]
\newtheorem{lem}[thm]{Lemma}
\newtheorem{prop}[thm]{Proposition}
\newtheorem{cor}[thm]{Corollary}

\usepackage{comment}
\usepackage{caption}

\begin{document}

\title{{\LARGE{Conley's index and connection matrices \\ for non-experts}}}

\author{
 \\
{~}\\
Phillipo Lappicy*\\
\vspace{2cm}}

\date{ }
\maketitle
\thispagestyle{empty}

\vfill

$\ast$\\
Instituto de Ciências Matemáticas e de Computação\\
Universidade de S\~ao Paulo\\
Avenida trabalhador são-carlense 400\\
13566-590, São Carlos, SP, Brazil\\


\newpage
\pagestyle{plain}
\pagenumbering{arabic}
\setcounter{page}{1}

\begin{abstract}
This is a self-contained tour of the Conley index and connection matrices. The starting point is Conley's fundamental theorem of dynamical systems. There is a short stop at the necessary topological background, before we proceed to the basic properties of the index. Then, the itinerary passes through the construction of connection matrices with a panoramic view of the applications: detect heteroclinic orbits arising in delay differential equations, and partial differential equations of parabolic type. The ride will be filled with examples and figures.

\ 

\end{abstract}

\section{Introduction}\label{sec:intro}

\numberwithin{equation}{section}
\numberwithin{figure}{section}
\numberwithin{table}{section}

Dynamical systems have played a huge role in modelling real life problems. From the motion of astronomical objects to the chemical reactions in our brain: dynamics is everywhere. A huge effort has been made in understanding what are typical dynamical systems, along with constructing all possible future asymptotic dynamics for a given system, i.e., construct the global attractor. See the survey of Bonatti \cite{Bonatti11} and references therein towards a global view of dynamical systems. 

For us, a dynamical systems is called a \emph{semiflow} $T(t)$ on the metric space $(X,d)$. Namely, a continuous map
\begin{align}\label{semiflow}
    T:\mathbb{R}_+\times X &\to X\\ \nonumber 
    (t,u_0)&\mapsto T(t)u_0
\end{align}
such that $T(0)u_0=u_0$ and $T(t+s)u_0=T(t)T(s)u_0$ for any $t,s\in \mathbb{R_+}$ and $u_0\in X$. 

Suppose that the semiflow $T(t)$ is \emph{bounded and dissipative}, i.e. any curve $T(t)u_0$ remains bounded for all times, and eventually enters some fixed large ball in $X$. Suppose also that \emph{orbits} given by $\{T(t)u_0 $ $ | $  $ t\in\mathbb{R}_+\}$ are \emph{precompact}, that is, the closure of orbits are compact in $X$. Hence, there exists the \emph{global attractor} $\mathcal{A}\subseteq X$, which is a compact set that attracts all bounded sets in $X$. In particular, $d(T(t)u_0,\mathcal{A})\rightarrow 0$ as $t\to\infty$ for all initial data $u_0 \in X$. See \cite{Hale88}, \cite{Lady91} and \cite{CaLaRo13}.

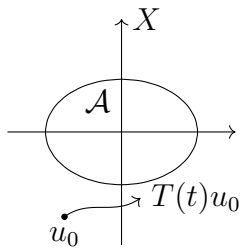
\begin{figure}[ht!]
\centering
\begin{tikzpicture}[scale=1]
    \draw[->] (-1.5,0) -- (1.5,0);
    \draw[->] (0,-1.5) -- (0,1.5) node[right] {$X$};

    \draw [domain=0:6.28,variable=\t,smooth] plot ({sin(\t r)},{0.7*cos(\t r)}) node[anchor= north east] {$\mathcal{A}$}; 
    
    \filldraw [black] (-0.75,-1.125) circle (1pt) node[anchor=north]{$u_0$};

    \draw [->,domain=-0.5:0.5,variable=\t,smooth] plot ({\t-0.25},{(\t)^3-1}) node[right] {$T(t)u_0$}; 
\end{tikzpicture}
\caption{The semiflow $T(t)$ is a time $t\in\mathbb{R}_+$ action of the space $X$, such that orbits of the action are curves in $X$ parametrized by time, displaying the time evolution of an initial point $u_0\in X$. For the bounded, dissipative and precompact case, any initial data converges to the global attractor.}
\end{figure}

Note that the semiflow is only well defined for positive times $t\in\mathbb{R}_+$, that is, not necessarily solutions are well defined backwards in time $t\in\mathbb{R}_-$, usually for the lack of uniqueness of backward solutions. 
Nevertheless, solutions are unique within the global attractor. In this case, the restriction of $T(t)$ to the compact subspace $\mathcal{A}$ is a \emph{flow}, and we can consider $t\in\mathbb{R}$. See \cite{CaLaRo13}.

In order to construct the (a priori quite complicated) global attractor, Morse's theory \cite{Morse34}, Smale's spectral decomposition \cite{Smale67} and Conley's fundamental theorem \cite{Conley78} were brilliant ideas to decompose the attractor into ``smaller sets", so that the problem of understanding the future asymptotics is now simpler: describe such smaller sets, and how they are related. In this decomposition, if there is complicated dynamics, it is confined to one of the smaller sets.


Before we present Conley's result, denote by $\mathcal{R}$ the \emph{chain-recurrent set}, which consists of all points $r_0\in \mathcal{A}$ such that for any $\epsilon>0$, there is a finite number of points $\{r_i\}_{i=0}^{n+1}\subseteq X$ and times $\{t_i\}_{i=0}^{n}\subseteq \mathbb{R}_+$ with $d(T(t_i)r_i,r_{i+1})<\epsilon$ for $i=0,...,n$ with $r_{n+1}:=r_0$.


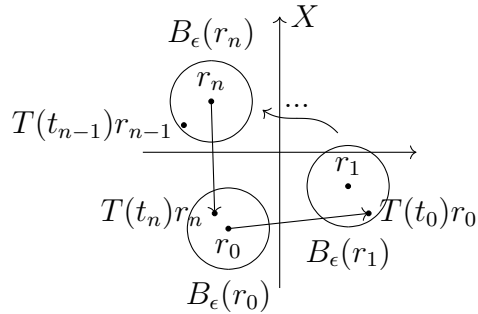
\begin{figure}[ht!]
\centering
\begin{tikzpicture}[scale=0.9]
    \draw[->] (-2,0) -- (2,0);
    \draw[->] (0,-2) -- (0,2) node[right] {$X$};

    \draw [domain=3.14:9.42,variable=\t,smooth] plot ({0.6*sin(\t r)-0.75},{0.6*cos(\t r)-1.125}) node[anchor= north] {$B_\epsilon (r_0)$}; 
    \filldraw [black] (-0.75,-1.125) circle (1pt) node[anchor=north]{$r_0$};
    \filldraw [black] (-0.95,-0.9) circle (1pt) node[anchor=east]{$T(t_{n})r_{n}$};        
    \draw[->] (-0.75,-1.125) -- (1.3,-0.9);

    \draw [domain=3.14:9.42,variable=\t,smooth] plot ({0.6*sin(\t r)+1},{0.6*cos(\t r)-0.5}) node[anchor= north] {$B_\epsilon (r_1)$}; 
    \filldraw [black] (1,-0.5) circle (1pt) node[anchor=south]{$r_{1}$};
    \filldraw [black] (1.3,-0.9) circle (1pt) node[anchor= west]{$T(t_0)r_{0}$};    

    \draw [domain=0:6.28,variable=\t,smooth] plot ({0.6*sin(\t r)-1},{0.6*cos(\t r)+0.75}) node[anchor= south] {$B_\epsilon (r_n)$}; 
    \filldraw [black] (-1,0.75) circle (1pt) node[anchor=south]{$r_{n}$};
    \filldraw [black] (-1.4,0.4) circle (1pt) node[anchor=east]{$T(t_{n-1})r_{n-1}$};    
    \draw[->] (-1,0.75) -- (-0.95,-0.9);    
    
    \draw [<-,domain=-0.5:0.6,variable=\t,smooth] plot ({\t+0.25},{-(\t)^3+0.5}); 
    \filldraw [black] (0.25,0.5) circle (0.01pt) node[anchor=south]{$...$};    
\end{tikzpicture}
\caption{A point $r_0\in \mathcal{R}$ eventually returns to a neighborhood of itself after some time of the flow has elapsed, where $B_\epsilon (r_i)$ denotes the ball in $X$ centered at $r_i$ and radius $\epsilon>0$.}
\end{figure}

\begin{thm}\emph{\textbf{Conley's fundamental theorem of dynamical systems.}}
Consider the semiflow $T(t)$ restricted to the compact attractor $\mathcal{A}\subseteq X$. Then, there is a continuous \emph{Lyapunov function} $L:\mathcal{A}\to \mathbb{R}$ which is strictly decreasing on $\mathcal{A} \backslash \mathcal{R}$. 
\end{thm}

In other words, any continuous semiflow on its compact global attractor  $\mathcal{A}$ is gradient-like off its chain-recurrent set $\mathcal{R}$. That is, all dynamics in $\mathcal{A} \backslash \mathcal{R}$ are connection orbits between connected components of $\mathcal{R}$, which are the ``smaller sets" mentioned above. 

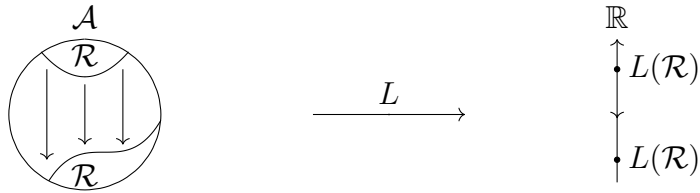
\begin{figure}[ht]\centering
\begin{tikzpicture}[scale=1]
    
    \draw [domain=0:6.28,variable=\t,smooth] plot ({sin(\t r)-7},{cos(\t r)}) node[anchor=south]{$\mathcal{A}$}; 
    \draw [domain=-0.575:0.575,variable=\t,smooth] plot ({\t-7},{(\t)^2+0.5});
    \filldraw [black] (-7,0.5) circle (0.01pt) node[anchor=south]{$\mathcal{R}$};
    \draw [domain=-0.725:0.75,variable=\t,smooth] plot ({\t-6.75},{(\t)^3-0.5});
    \filldraw [black] (-7,-0.515) circle (0.01pt) node[anchor=north]{$\mathcal{R}$};
    
    \draw[->] (-7,0.4) -- (-7,-0.4);
    \draw[->] (-6.5,0.6) -- (-6.5,-0.4);
    \draw[->] (-7.5,0.6) -- (-7.5,-0.6);
            
    \draw[->] (-4,0) -- (-2,0);
    \filldraw [black] (-3,0) circle (0.1pt) node[anchor=south]{$L$};

    \draw[->] (0,-0.9) -- (0,1) node[above] {$\mathbb{R}$};
    
    \filldraw [black] (0,0.6) circle (1pt) node[anchor= west]{$L(\mathcal{R})$};    
    \draw[->] (0,0) -- (0,-0.05) ;
    \filldraw [black] (0,-0.6) circle (1pt) node[anchor= west]{$L(\mathcal{R})$};

\end{tikzpicture}
\captionof{figure}{Dynamcal decomposition of the global attractor $\mathcal{A}$: gradient structure off the chain-recurrent set $\mathcal{R}$, induced by the one-dimensional dynamics arising from the Lyapunov function $L$.}
\end{figure}

Then, in order to understand the global attractor $\mathcal{A}$, we need to describe the set $\mathcal{R}$, and know which connected components of the recurrent set are connected through a connection orbit in $\mathcal{A}\backslash \mathcal{R}$. 

For now, suppose that the recurrent set $\mathcal{R}$ can be decomposed into $n$ disjoint ``energy levels", namely $\mathcal{R}=\cup_{i=1}^n M_i$. Moreover, say that the $i^{th}$ energy level has $C_i$ connected components, $M_i=\cup_{i'=1}^{C_i} M_{i,i'}$. A tool which is capable of detecting interesting dynamical behaviour in such smaller sets $M_{i,i'}$ within the global attractor $\mathcal{A}$ is now called the \emph{Conley index}. The index associates a homotopy type to the set $M_{i,i'}$ that relates the dynamics of a neighborhood of $M_{i,i'}$ to the dynamics of its neighborhood boundary, displaying its nearby unstable dynamics. This will be further discussed in Section \ref{sec:index}, where we compute the index for some possibilities that lie in the chain recurrent set: hyperbolic fixed points, periodic orbits, homoclinic orbits, and other examples. 


The question that remains is: what can we say about $\mathcal{A} \backslash \mathcal{R}$? That is, is there a \emph{heteroclinic connection} $T(t)u_0$ between $M_{i,i'}$ and $M_{j,j'}$ given by
\begin{equation}
M_{i,i'} \xleftarrow{t\to-\infty} T(t)u_0 \xrightarrow{t\to\infty} M_{j,j'}?
\end{equation}

In order to answer this question, we construct the \emph{connection matrix} $\Delta$ such that its entries $\Delta_{i,j}$ are also matrices that relate the Conley index of $M_i$ and $M_j$. 
As the below theorem indicates, each entry $\delta_{i',j'}$ of the matrix $\Delta_{i,j}$ encodes the information of the existence of connection orbits from $M_{i,i'}$ to $M_{j,j'}$ and consequently, the structure of the set $\mathcal{A} \backslash \mathcal{R}$. Moreover, knowledge of the Conley index of the smaller sets $M_{i,i'}$ and properties of the connection matrix $\Delta$ are enough to reconstruct most its entries. Those will be explored in Section \ref{sec:matrix}.

\begin{thm}\emph{\textbf{Heteroclinic detection.}}\label{thmdetectionintro}
If the entry $\delta_{i^{'},j^{'}}$ of the matrix $\Delta_{i,i+1}$ is nonzero, then there is a heteroclinic orbit from $M_{i+1,j'}$ to $M_{i,i'}$, i.e., a connection from the $j^{'}$ connected component of $M_{i+1}$ to the $i^{'}$ connected component of $M_{i}$.
\end{thm}



There are three applications that will be treated in these notes: the heteroclinics between equilibria and/or periodic orbits of given unstable dimension, the connections of the Chafee-Infante attractor, and the general connection problem for partial differential parabolic equations in one spatial dimension. See Section \ref{sec:app}.

\begin{figure}[ht]
\minipage{0.5\textwidth}\centering%
    \begin{tikzpicture}[scale=1]

    \filldraw [black] (0,4) circle (2pt) node[anchor=south] {$u\equiv 0$}; 
    
    \draw[thick,->] (0,4) -- (0,3.1);

    \draw (-1,2.5) arc (180:360:1cm and 0.5cm);
    \draw (-1,2.5) arc (180:0:1cm and 0.5cm);
    \draw[->] (-1,2.5) -- (-1,2.49) node[anchor=east]{$M_1$};
    
    \draw[thick,->] (0,2) -- (0,1.1);

    \draw (-1,0.5) arc (180:360:1cm and 0.5cm);
    \draw (-1,0.5) arc (180:0:1cm and 0.5cm);
    \draw[->] (-1,0.5) -- (-1,0.49) node[anchor=east]{$M_0$};
    
    \end{tikzpicture}
\endminipage\hfill
\minipage{0.5\textwidth}\centering%
\begin{tikzpicture}[scale=0.8]
\filldraw [black] (0,0) circle (3pt) node[anchor=south]{$e_{N}$};
\filldraw [black] (-1,-1) circle (3pt) node[anchor=east]{$e^+_{n}$};
\filldraw [black] (1,-1) circle (3pt) node[anchor=west]{$e^-_{n}$};

\filldraw [black] (-1,-2.4) circle (0.5pt);
\filldraw [black] (-1,-2.5) circle (0.5pt);
\filldraw [black] (-1,-2.6) circle (0.5pt);
\filldraw [black] (1,-2.4) circle (0.5pt);
\filldraw [black] (1,-2.5) circle (0.5pt);
\filldraw [black] (1,-2.6) circle (0.5pt);

\filldraw [black] (-1,-4) circle (3pt) node[anchor=east]{$e^+_{1}$};
\filldraw [black] (1,-4) circle (3pt) node[anchor=west]{$e^-_{1}$};
\filldraw [black] (-1,-5) circle (3pt) node[anchor=east]{$e^+_{0}$};
\filldraw [black] (1,-5) circle (3pt) node[anchor=west]{$e^-_{0}$};

\draw[thick,->] (0,0) -- (-0.9,-0.9);
\draw[thick,->] (0,0) -- (0.9,-0.9);

\draw[thick,->] (-1,-1) -- (0.9,-1.9);
\draw[thick,->] (1,-1) -- (-0.9,-1.9);
\draw[thick,->] (-1,-1) -- (-1,-1.87);
\draw[thick,->] (1,-1) -- (1,-1.87);

\draw[thick,->] (-1,-3) -- (0.9,-3.9);
\draw[thick,->] (1,-3) -- (-0.9,-3.9);
\draw[thick,->] (-1,-3) -- (-1,-3.87);
\draw[thick,->] (1,-3) -- (1,-3.87);

\draw[thick,->] (-1,-4) -- (0.9,-4.9);
\draw[thick,->] (1,-4) -- (-0.9,-4.9);
\draw[thick,->] (-1,-4) -- (-1,-4.87);
\draw[thick,->] (1,-4) -- (1,-4.87);

\end{tikzpicture}
\endminipage
\caption{On the left, the global attractor $\mathcal{A}$ for a delay differential equation has heteroclinics between the equilibrium $u\equiv 0$ and the periodic orbit $M_1$, and also between the periodic orbits $M_1$ and $M_0$. On the right, the global attractor $\mathcal{A}$ for the Chafee-Infante equation has heteroclinics between $e_N$ and $e^\pm_{n}$, and also between $e^\pm_{k}$ and $e^\pm_{k-1}$ for any $k=1,...,n$.}
\label{fig:attINTRO}
\end{figure}
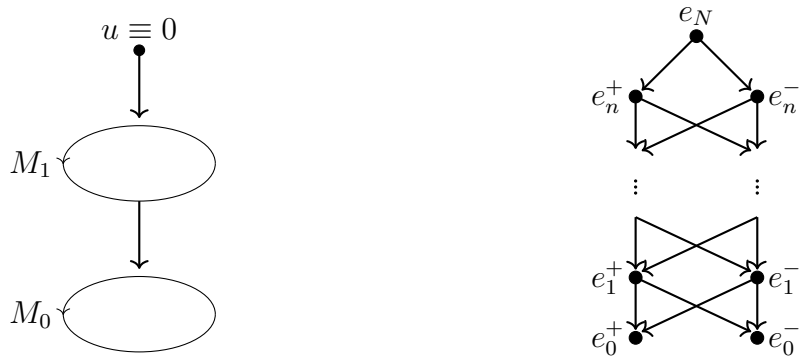

We mention a brief historical timeline regarding the results of this survey. After Conley's general framework was originally set with Easton in \cite{ConleyEaston71}, Rybakowski extended the theory for infinite dimensional dynamical systems in \cite{Rybakowski82} and \cite{Rybakowski87}, then Franzosa proved the existence of a connection matrix \cite{Franzosa86} and \cite{Franzosa89} in order to detect heteroclinics. For a deeper yet brief historical account, see \cite{Mischaikow99}.

The main references of this paper are Conley's seminal exposition \cite{Conley78}, Rybakowski \cite{Rybakowski82}, Smoller \cite{Smoller83}, Mischaikow \cite{Mischaikow87}, \cite{Mischaikow99}, \cite{MischaikowMrozek02}. 

We do not want to survey the whole theory of Conley and its ramifications. We leave out several realms still to be explored outside these notes, such as transition matrices \cite{FranzosaMischaikow98} and Morse Inequalities \cite{Smale60}. The former can be used in order to understand the structure of global bifurcations, and changes in the connection matrix, as the unified theory of Franzosa and Vieira \cite{FranzosaVieira17}. The latter was adapted by Conley and Zehnder \cite{ConleyZehnder83} in order to prove Arnold's conjecture that counts the number of fixed point of measure preserving diffeomorphisms of the torus, and generalized by Floer \cite{Floer88} to prove Arnold's conjecture for symplectic maps \cite{Floer86}. It was also used (when the inequality is indeed an equality) to show Poincaré's conjecture for dimension bigger than 5 by Smale \cite{Smale61}. 

Instead, we focus only on the index and the connection matrices. We seek to complement the existing material, collecting known examples, and explaining the geometric ideas while omitting most proofs. Nevertheless, we present the akin relationship between the algebraic proof of the zig-zag lemma from topology, and its geometric counterpart with dynamical interpretation that arises in the construction of the connection matrix.

The remaining of this survey is organized as follows. In Section \ref{sec:alg} we introduce the basics of algebraic topology, including exact sequences, CW-complexes and cellular homology theory. In Section \ref{sec:index} we develop the Conley index along with its properties, and compute simple examples. In Section \ref{sec:matrix} we describe the connection matrix, its properties, and its construction in a particular case. Lastly, in Section \ref{sec:app}, we apply the properties of the matrix in order to obtain heteroclinic connections between invariant sets for delay and parabolic differential equations.

\textbf{Acknowledgment.} I am thankful for Alexandre N. Carvalho and Bernold Fiedler that have encouraged diving into this topic, and taught me the ``basics" so that I had the necessary equipment for the nosedive. Maria C. Carbinatto helped with several discussions and suggestions that followed the original notes from Konstantin Mischaikow's courses in the 90's, that have several examples incorporated here. The participants from the minicourse taught at ICMC Summer School 2019 that had this note as an outcome had fundamental remarks. In particular, Alexandre do Nascimento, Sérgio Monari, and Tiago Pereira. The author was supported by FAPESP, Brasil, grant number 2017/07882-0.

\section{Algebraic topology at a glance}\label{sec:alg}

A common language used between topologists is of exact sequences and chain complexes, that relate spaces through morphisms (generalization of maps) in diagram form. In particular, some topological spaces can be described by CW-complexes, which consists of building blocks called cells, and recipes on how those cells are topologically glued together in order to form the topological space. From these, we obtain the cellular homology theory, an algebraic structure that accounts for topological properties of the space. Therefore, we proceed to introduce some notation about exact sequences, and cellular homology. This section was mostly excised out of \cite{Dieck08}.

A \emph{sequence} is a schematic way to pin down the relationship between spaces and morphisms altogether through a diagram that concatenates the objects $A_n$ and their morphisms $a_n:A_n\to A_{n-1}$ as
\begin{equation}\label{Seqdef}
    ...\rightarrow A_n \xrightarrow{a_n} A_{n-1} \xrightarrow{a_{n-1}} A_{n-2} \rightarrow ... \hspace{0.2cm} .
\end{equation}

We say the sequence \eqref{Seqdef} is \emph{exact} if the image of each morphism is equal to the kernel of the next: $Im(a_n)=Ker(a_{n-1})$. In particular, the sequence $0\rightarrow A_{n}\rightarrow_{a_{n}} A_{n-1}$ is exact, if and only if, $a_{n}$ is injective. We denote by $0$ the trivial object, let it be a group, ring, module or others. Similarly, the sequence $A_{n}\rightarrow_{a_{n}} A_{n-1}\rightarrow 0$ is exact, if and only if, $a_{n}$ is surjective. Concatenating both these sequences yields a schematic way to say the morphism $a_n$ is bijective.

A \emph{short exact sequence} is a given exact sequence with three nontrivial elements: 
\begin{equation}\label{SESdef}
    0\rightarrow A_n \xrightarrow{a_n} A_{n-1} \xrightarrow{a_{n-1}} A_{n-2} \rightarrow0.
\end{equation}

This means $a_{n}$ is injective, $a_{n-1}$ is surjective, and $Im (a_{n}) = Ker (a_{n-1})$. It is useful to think of the short exact sequence \eqref{SESdef} as way to decompose the space $A_{n-2}$ as the quotient of $A_{n-1}$ with respect to the image of $A_n$ through $a_n$ given by the following isomorphism 
\begin{equation}
    A_{n-2}\cong Im (a_{n-1}) \cong \frac{A_{n-1}}{Ker (a_{n-1})}\cong \frac{A_{n-1}}{Im (a_{n})} .
\end{equation}

An example is the inclusion of the circle in the boundary of a disk, and identifying the boundary of the disk into one point in order to obtain a two-sphere, yielding the following short exact sequence
\begin{equation}\label{SESsphere}
    0\rightarrow \mathbb{S}^1 \xrightarrow{inc} B^1 \xrightarrow{proj} \frac{B^1}{\mathbb{S}^1}\cong \mathbb{S}^2 \rightarrow 0
\end{equation}
where $inc$ is injective, $proj$ is surjective, and $Im(inc)=Ker (proj)$.

We say the short exact sequence \eqref{SESdef} \emph{splits} if there exists either a map $d:A_{n-1}\to A_n$ such that $d \circ a_n =id_{A_n}$, or $\tilde{d}:A_{n-2}\to A_{n-1}$ such that $a_{n-2}\circ \tilde{d} =id_{A_{n-2}}$. In case a short exact sequence splits, then the splitting lemma says that $A_{n-1}=A_n\oplus A_{n-2}$.

An example of a split short exact sequence is the inclusion of a circle in the wedge sum of said circle and a sphere, and then projecting into the sphere:
\begin{equation}
    0\rightarrow \mathbb{S}^1 \xrightarrow{inc} \mathbb{S}^1\vee \mathbb{S}^2 \xrightarrow{proj} \mathbb{S}^2 \rightarrow 0
\end{equation}
where $inc$ is injective, $proj$ is surjective, and $Im(inc)=Ker (proj)$. 

We recall the \emph{wedge sum} is a ``one-point gluing" of topological spaces: if $Y$ and $Z$ are topological spaces with  $y_0\in Y$ and $z_0\in Z$, the wedge sum $Y\vee Z$ is the quotient space of the disjoint union of $Y$ and $Z$ by the identification $y_0\sim z_0$.

A \emph{long exact sequence} is an exact sequence with more than three nontrivial elements: 
\begin{equation}\label{LESdef}
    ... \rightarrow A_n \xrightarrow{a_n}A_{n-1} \xrightarrow{a_{n-1}}A_{n-2} \rightarrow ...
\end{equation}

Next, we introduce the geometric notion of CW-complexes.

An \emph{attaching map} is a continuous map $f:\partial B^n\to Y$ that will attach an $n$-dimensional ball $B^n$ to the topological space $Y$ by relating its boundary $\partial B^n$ and $f(\partial B^n)$ as follows. The \emph{adjunction space} $B^n\cup_{f} Y:= B^n\sqcup Y / \sim $, where $\sqcup$ denotes the disjoint union of sets, and $\sim$ is an equivalence relation that idenfities $\partial B^n$ and $f(\partial B^n)$ pointwise: given $b\in \partial B^n,y\in Y$, then $b\sim y$ if, and only if $f(b)=y$. The resulting $c^n:=(B^n\cup_{f} Y)\backslash Y$ is called an \emph{$n$-dimensional (open) cell}.  
If $f$ is a homeomorphism, and the closure of each $n$-dimensional open cell is homeomorphic to $B^n$, then we say the attaching map is \emph{regular}.

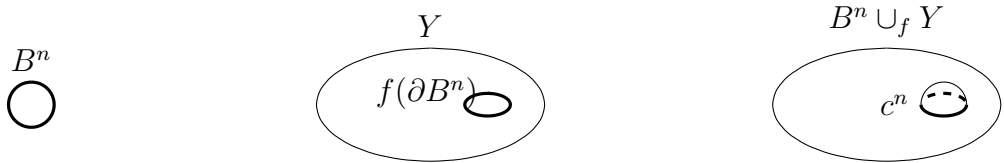
\begin{figure}[ht]\centering
    \begin{tikzpicture}[scale=1.5]

        \draw [very thick,domain=0:6.28,variable=\t,smooth] plot ({-2.5+0.2*sin(\t r)},{0.2*cos(\t r)});
        \draw [domain=0:6.28,variable=\t,smooth] plot ({-2.5+0.2*sin(\t r)},{0.2*cos(\t r)}) node[anchor= south] {$B^n$};        

        \draw [domain=0:6.28,variable=\t,smooth] plot ({1+sin(\t r)},{0.5*cos(\t r)}) node[anchor= south] {$Y$};
        
        \draw [very thick,domain=0:6.4,variable=\t,smooth] plot ({1.5+0.2*sin(\t r)},{0.1*cos(\t r)}) node[anchor= east] {$f(\partial B^n)$};
        
        \draw [domain=0:6.28,variable=\t,smooth] plot ({5+sin(\t r)},{0.5*cos(\t r)}) node[anchor= south] {$B^n\cup_f Y$};
        
        \draw [very thick,domain=1.57:4.71,variable=\t,smooth] plot ({5.5+0.2*sin(\t r)},{0.1*cos(\t r)}) node[anchor= east] {$c^n$};
        \draw [very thick, dashed,domain=5.5:7.65, variable=\t,smooth] plot ({5.5+0.2*sin(\t r)},{0.1*cos(\t r)}) ;
        \draw [domain=4.71:7.85,variable=\t,smooth] plot ({5.5+0.2*sin(\t r)},{0.2*cos(\t r)});

    \end{tikzpicture}
\caption{We glue the $n$-dimensional ball $B^n$ through its boundary $\partial B^n$ to the set $f(\partial B^n)$ within the topological space $Y$, yielding the $n$-dimensional cell $c^n$ in $B^n\cup_f Y$.}
\end{figure}

An example, we consider the topological space which is given by a point $Y=\{ c^0\}$. To this $0$-dimensional cell, we glue an interval through the attaching map $f:\partial [0,1]\to Y$ such that $f(0)=f(1)=c^0$. Hence, $0$ and $1$ lie in the same equivalence class of $c^0$, and the adjunction space $[0,1]\cup_f Y$ is given by a circle $\mathbb{S}^1$ with distinguished point $c^0$, and a $1$-dimensional cell $c^1:=\mathbb{S}^1\backslash Y$.

Inductively, we define a \emph{CW-complex} of a topological space $Y$ of dimension $n$. Begin with $N_0$ cells of dimension $0$, given by $Y_0:=\{c^0_j\}_{j=1}^{N_0}$, and glue $N_1$ cells of dimension $1$ to $Y_0$ in order to obtain $Y_1:=(\sqcup_{j=1}^{N_1} B^1)\cup_{f_1} Y_0$. We continue this process inductively until we glue $N_n$ cells of dimension $n$ to $Y_{n-1}$, to obtain $Y_n:=(\sqcup_{j=1}^{N_n} B^n)\cup_{f_n} Y_{n-1}$. This yields the nested sequence
\begin{equation}\label{filtX}
    Y_{-1} \subseteq  Y_0\subseteq Y_1\subseteq ... \subseteq Y_n=Y
\end{equation}
where $Y_{-1}=\varnothing$ for well definition purposes, and $Y_k$ is called the \emph{$k^{th}$-skeleton of $Y$}. Moreover, if all attaching maps are regular, then the CW-complex $Y$ is \emph{regular}. 

Few remarks are at hand. Firstly the quotient $Y_k / Y_{k-1}$ consists of the wedge of all $k$-dimensional cells $\vee_{j=1}^{N_k} c^k_j$, with distinguished point being the equivalence class of $Y_{k-1}$, where $N_k$ denotes the total number of $k$-dimensional cells. Secondly, we allow a CW-complex to have no $k$-cell at all, for any $k\geq 0$. 
Thirdly, there is no unique way to decompose a space $Y$ into a CW-complex. 

As an example, we provide two different CW-complex decompositions of the sphere.

The first is obtained by starting with a point $Y_{0}=\{c^0\}$, then adding no $1$-cells, and hence $Y_{1}=Y_{0}$. Subsequently one $2$-cell $\{c^2\}$ is glued to $Y_1$, where the attaching map $f:\partial B^2\to Y_1$ is given by the constant map $f(\partial B^2)=c^0$, which geometrically is the collapse of the circle as the boundary of a disk into the point $c^0$.

The second is obtained by starting with no zero dimensional cells $Y_{0}=\varnothing$, but one circle denoting the equator as the $1$-cell $Y_{1}=\{c^1\}$. This is given by the attaching map $f_1: \partial [0,1]\to \partial [0,1]$ such that $f_1(0)=f_1(1)$. Subsequently, we glue two $2$-cells $\{c^2_+,c^2_-\}$ denoting the upper and lower hemispheres. The attaching map is given by the identification of the circle as the boundary of both disks $c^2_\pm$ with the equator $c^1$, through $f:\partial B^2\to Y_1$ such that $f(\partial B^2)=c^1$, the identity map on the circle. 
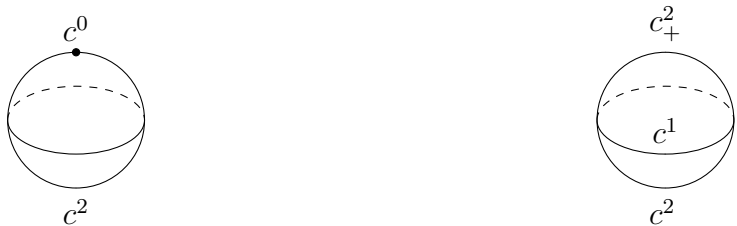
\begin{figure}[ht]
\minipage{0.5\textwidth}\centering%
    \begin{tikzpicture}[scale=0.9]

        \draw (-1,0) arc (180:360:1cm and 0.5cm);
        \draw[dashed] (-1,0) arc (180:0:1cm and 0.5cm);
        
        \draw (0,0) circle (1cm);
    
        \filldraw [black] (0,1) circle (1.5pt) node[anchor=south] {$c^0$}; 
        \filldraw [black] (0,-1) circle (0.01pt) node[anchor=north] {$c^2$};     
    \end{tikzpicture}
\endminipage\hfill
\minipage{0.5\textwidth}\centering%
    \begin{tikzpicture}[scale=0.9]

        \draw (-1,0) arc (180:360:1cm and 0.5cm);
        \draw[dashed] (-1,0) arc (180:0:1cm and 0.5cm);
        \draw[-] (-0.01,-0.5) -- (0,-0.5) node[anchor=south] {$c^1$};
        
        \draw (0,0) circle (1cm);
    
        \filldraw [black] (0,1) circle (0.01pt) node[anchor=south] {$c^2_+$}; 
        \filldraw [black] (0,-1) circle (0.01pt) node[anchor=north] {$c^2_-$};
    \end{tikzpicture}
\endminipage
\caption{Two CW-complex structures of the sphere: on the left, a point $c^0$ being the north pole with a 2-cell $c^2$ glued to it; on the right, a circle $c^1$ being the equator with two 2-cells $c^2_+$ and $c^2_-$ being the hemispheres glued to it.}
\end{figure}

Even though this hierarchical description of a topological space is widely spread in by topologists, we mention a dynamical systems perspective of them. In the particular case the topological space $Y$ is the global attractor $\mathcal{A}$ of a parabolic equation on the interval, we know that $\mathcal{A}=\cup_{i=1}^n W^u_k(e_i)$ where $e_i$ are the hyperbolic equilibria, and $W^u_k(e_i)$ its corresponding unstable manifold with dimension $k$. It was proved, by Fiedler and Rocha \cite{FiedlerRocha15}, that actually $\mathcal{A}$ is a regular CW-complex such that its $k$-cells are $c^k_i:=W^u_k(e_i)$ and 0-cels are the equilibria $\{e_i\}_{i=1}^n$. See also the survey of Bott \cite{Bott88}. Therefore, certain complicated topological spaces which do not have a CW-structure can not occur as global attractors of those dynamics. 

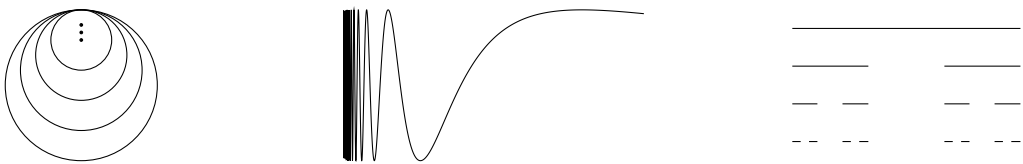
\begin{figure}[ht]
\minipage{0.3\textwidth}\centering%
    \begin{tikzpicture}[scale=1]
        \draw (0,0) circle (1cm);
        \draw (0,0.2) circle (0.8cm);
        \draw (0,0.4) circle (0.6cm);
        \draw (0,0.6) circle (0.4cm);
        \filldraw [black] (0,0.8) circle (0.5pt); 
        \filldraw [black] (0,0.7) circle (0.5pt); 
        \filldraw [black] (0,0.6) circle (0.5pt); 
    \end{tikzpicture}
\endminipage\hfill
\minipage{0.3\textwidth}\centering%
    \begin{tikzpicture}[scale=1]
        \draw[domain=0.01:0.8,samples=5000] plot (5*\x, {sin((1/\x)r)});
    \end{tikzpicture}
\endminipage\hfill
\minipage{0.3\textwidth}\centering%
    \begin{tikzpicture}[scale=1]
        \draw (0,0) -- (3,0);
        
        \draw (0,-0.5) -- (1,-0.5);        
        \draw (2,-0.5) -- (3,-0.5); 
        
        \draw (0,-1) -- (0.33,-1);
        \draw (0.66,-1) -- (1,-1);
        \draw (2,-1) -- (2.33,-1);
        \draw (2.66,-1) -- (3,-1);

        \draw (0,-1.5) -- (0.11,-1.5);
        \draw (0.22,-1.5) -- (0.33,-1.5);
        \draw (0.66,-1.5) -- (0.77,-1.5);
        \draw (0.88,-1.5) -- (1,-1.5);

        \draw (2,-1.5) -- (2.11,-1.5);
        \draw (2.22,-1.5) -- (2.33,-1.5);
        \draw (2.66,-1.5) -- (2.77,-1.5);
        \draw (2.88,-1.5) -- (3,-1.5);        
    \end{tikzpicture}
\endminipage
\caption{On the left, the Hawaiian earring, on the middle the topologist sine curve, on the right a construction of the Cantor set: none of those can be attractors of parabolic equations.}\label{fig}
\end{figure}

This raises the question: for which evolution equations with global attractors can one obtain a CW-complex decomposition of $\mathcal{A}$? 
Moreover, what are the restrictions the semiflow impose on the topological structures of the attractor through its CW-complex structure? Which kind of evolution equations yield more intricate topological objects as attractors, such as the Möbius band, the Klein bottle, Alexander's horned sphere or one of the spaces in Figure \ref{fig}? See \cite{SanchezGabites11} for a partial answer in dimension 3. But we digress with those interesting open questions...

Nevertheless, we mention that the filtration \eqref{filtX} already gives hints that a topological tool should yield information about the connection between sets in phase-space of a dynamical systems. Indeed, let $Y$ be the global attractor $\mathcal{A}$ with filtration \eqref{filtX}, where the $0$-skeleton are the equilibria points $\{e_i\}_{i=1}^n$. Consider two equilibria $e_i$ and $e_j$, and define the space $\Sigma:= \{e_i,e_j\}\cup (W^u(e_i)\cap W^s(e_j))$ as a subcomplex of $\mathcal{A}$. Then, a connection between two equilibria is reflected in the quotient $\Sigma / \{e_i,e_j\}$. That is, if there is no connection, then the quotient has two connected components, whereas if there is a connection, then the same quotient has one connected component. Moreover, the fact that $W^u(e_i)\cap W^s(e_j)$ is $\varnothing$ or not will be reflected in the creation of a cycle of $\Sigma / \{e_i,e_j\}$, or not. That is, this information can be captured by homology, which is the theory we develop next. 

We now introduce the cellular complex of a CW-complex, which uses the algebraic language of sequences in order to explore its geometric notion. Associated to the chain complex, there is its homology group.

Given a CW-complex structure of $Y$, consider the vector space $C_k(Y)$ with coefficients in $\mathbb{Z}_2$, where its basis is given by the $k$-dimensional cells $\{c^k_j\}_{j=1}^{N_k}$ of $Y$,
\begin{equation}\label{cellcomplex}
    C_k(Y):=\left\{\sum_{j=1}^{N_k} z_j c^k_j \text{ $|$ } z_j\in \mathbb{Z}_2\right\} \cong \mathbb{Z}_2^{N_k}.
\end{equation}

We will denote by $\{0\}$ the trivial vector space from now on.

The \emph{cellular complex} $(C_k(Y),\partial_k)$ of the topological space $Y$ with dimension $n$ is given by the vector space $C_k(Y)$ as in \eqref{cellcomplex}, and homomorphisms $\partial_k: C_k(Y)\to C_{k-1}(Y)$ called \emph{boundary maps} with the property that $\partial_k \circ \partial_{k+1}=0$ for all $k$, or written in short as $\partial^2=0$. Note that the boundary maps $\partial_k$ decreases dimension by 1. Schematically, it can be written as the sequence  
\begin{equation}\label{CELLLES}
    C_n(Y) \xrightarrow{\partial_n}C_{n-1}(Y) \xrightarrow{\partial_{n-1}}... \xrightarrow{\partial_{1}} C_0(Y) \xrightarrow{\partial_{0}} 0.
\end{equation}

The formula for the boundary maps $\partial _k$ can computed from its knowledge on each basis element $c^k_i\in C_k(X)$, since they are homomorphisms, and is given by
\begin{equation}
    \partial_k (c^k_i)=\sum_{j=1}^{N_{k-1}} \text{deg}(f_j) c^{k-1}_j    
\end{equation}
where $\text{deg}(f_j)$ denotes the usual degree of the attaching map $f_j$, that glues $B^k$ to $c^{k-1}_j$.

Naturally, elements of $Im (\partial_{k+1})$ are called \emph{(k+1)-boundaries}, since the operator $\partial_{k+1}$ maps the $(k+1)$-dimensional cells to their respective $k$-dimensional boundaries. On the other hand, elements of $Ker (\partial_k)$ are called \emph{k-cycles}, since those are the $k$-dimensional cells that have their whole boundary identified to one point denoted by 0. Due to the condition $\partial_k \circ \partial_{k+1}=0$, all $(k+1)$-boundaries are $k$-cycles. That is, the image of $\partial_{k+1}$ is contained in the kernel of $\partial_{k}$.

Since the images are contained in the kernels, we can consider the quotient space that defines the group of $k$-dimensional cycles with their boundaries as equivalence classes, known as \emph{$k^{th}$-dimensional homology group},
\begin{equation}
    H_k(Y):=\frac{Ker (\partial_k)}{Im (\partial_{k+1})}.
\end{equation}

We give an algebraic and geometric interpretation of the homology groups $H_k(Y)$. Recall that the sequence \eqref{CELLLES} is exact if $Im(\partial_{k+1})=Ker(\partial_{k})$. Algebraically, the homology groups of $Y$ measure ``how far" the chain complex associated to $Y$ is from being exact. Geometrically, the $k^{th}$-homology computes how complicated the space $Y$ is in $k^{th}$-dimension. In particular, its elements constitute the $k^{th}$-dimensional ``holes" of the space $Y$. In other words, a nontrivial element of the $k^{th}$-homology is a $k$-cycle which is not the boundary of any $k+1$ subspace. A hypothetical $k+1$  space with such $k^{th}$-homology cycle as the boundary would be the ``space which not there", that is, the hole. In particular, the $0^{th}$-dimensional homology counts the connected components of the set $Y$.

We now compute the homology groups of the sphere, given a previously chosen CW-complex structure. We note that the same homology groups are obtained for the two different CW-decompositions. In general, it can be shown that the homology is independent of the CW-complex structure. 

Consider the following CW-structure of the sphere
\begin{equation}\label{CWsphere}
    0 \xrightarrow{0} C_2(Y) \xrightarrow{\partial_2} C_1(Y) \xrightarrow{\partial_1} C_0(Y) \xrightarrow{0} 0
\end{equation}
where $C_2(Y)=\text{span}\{c^2\}\cong \mathbb{Z}_2$, $C_1(Y)=0$ and $C_0(Y)=\text{span}\{c^0\}\cong \mathbb{Z}_2$. Therefore, $\partial_2\equiv 0$ and $\partial_1\equiv 0$. This implies that
\begin{align*}
    H_2(Y)&=\frac{Ker (\partial_2)}{Im (0)}=C_2(Y)\cong \mathbb{Z}_2\\
    H_1(Y)&=\frac{Ker (\partial_1)}{Im (\partial_{2})}= 0\\
    H_0(Y)&=\frac{Ker (0)}{Im (\partial_{1})}=C_0(Y)\cong \mathbb{Z}_2.
\end{align*}

In the above algebraic interpretation, the CW-structure of the sphere in the sequence \eqref{CWsphere} lacks to be exact in dimension 0 and 2. Geometrically, it has one connected component, and one ``hole" in dimension 2.


We mention few properties of the cellular homology. Recall that $Y_{k}/Y_{k-1}= \vee_{i=1}^{N_k} c^k_i$. Therefore its only nonzero homology is $H_k(Y_{k}/Y_{k-1})=\bigoplus_{i=1}^{N_k} \mathbb{Z}_2$. Also, higher dimensional cells do not contribute fo lower dimensional homology, i.e., $H_k(Y)\cong H_k(Y_j)$ for all $k<j$. Lastly, if $Y$ has no $k$-cells, 
then the $k^{th}$-homology is trivial $H_k(Y)=0$. In particular, if the dimension of $Y$ is $n$, then all $k^{th}$-homology with $k>n$ are trivial.

Below we mention the fundamental properties of homological theory. For some, these are treated as the Eilenberg–Steenrod axioms, and hence the starting point to study abstract nonsense homology theories.

\begin{thm}\label{thm:H_k}
Consider two CW-complexes $Y,Z$ with respective $j$-skeleton $Y_j,Z_j$ for all $j\in\mathbb{N}$. Then, the following properties hold. 
\begin{enumerate}
    \item \textbf{\emph{Homotopy invariance:}} if $Y$ and $Z$ are homotopy equivalent, then for all $k$
    \begin{equation*}
        H_k(Y)\cong H_k(Z).
    \end{equation*}
    
    \item \textbf{\emph{Excision:}}  let  $U\subseteq A\subseteq Y$ such that $A$ is closed, $U$ is open, $clos(U)\subseteq int(A)$, and $A$ is a deformation retract in $U$. Then, for all $k$, we have that
    \begin{equation*}
        H_k\left(\frac{Y}{A}\right)\cong H_k\left(\frac{Y\backslash U}{A\backslash U}\right).
    \end{equation*}
    
    \item \textbf{\emph{Dimension:}}  $H_{k}(Y_0)=0$ for all $ k\neq 0$.
    
    \item \textbf{\emph{Additivity:}}  If $Y = \coprod_{j}{Y_{j}}$, the disjoint union of topological spaces $Y_{j}$, then 
    \begin{equation*}
        H_k(Y) \cong \bigoplus_{j} H_k(Y_{j}).
    \end{equation*}
    
    \item \textbf{\emph{Exactness:}} The inclusion $i:Y_{j-1}\to Y_j$ and projection $p: Y_j \to Y_j/Y_{j-1}$ induce the long exact sequence
    
\begin{center}
\begin{tikzpicture}
\matrix[matrix of nodes,ampersand replacement=\&, column sep=0.5cm, row sep=0.5cm](m)
{
 \& $\cdots$ \& $H_{k+1} \left(\frac{Y_j}{Y_{j-1}}\right)$ \\
$H_k(Y_{j-1})$ \& $H_k(Y_j)$ \& $H_k \left(\frac{Y_j}{Y_{j-1}}\right)$ \\
$H_{k-1}(Y_{j-1})$ \& $\cdots$ \&  \\
};
\draw[->] (m-1-2) edge (m-1-3)
          (m-1-3) edge[curvedlink] (m-2-1)
          (m-2-1) edge (m-2-2) 
          (m-2-2) edge (m-2-3)
          (m-2-3) edge[curvedlink] (m-3-1)
          (m-3-1) edge (m-3-2);
\filldraw [black] (4.1075,-0.5) circle (0.01pt) node[anchor=west] {$d_{k-1}$};
\filldraw [black] (4.265,1) circle (0.01pt) node[anchor=west] {$d_{k}$};
\end{tikzpicture}.
\end{center}

\end{enumerate}
\end{thm}

Besides the cellular homology described above, there is a theory known as reduced homology, which agrees with the cellular homology in all dimensions, except in dimension 0. In dimension 0, the reduced homology equals the cellular homology quotiented by $\mathbb{Z}_2$. For application purposes, we will deal with reduced homology in the upcoming sections. We note that the remark in this paragraph is only true for regular pairs, which will be the case for the Conley index in the next sections.

As the exactness property in Theorem \ref{thm:H_k}, there is an interesting way to obtain a long exact sequence from a short exact sequences, as an application of the snake lemma. This will be the heart of the construction of the connection matrix.

\begin{lem}\emph{\textbf{Zig-Zag lemma.}}\label{thmzigzag}
Consider the short exact sequence
\begin{equation}\label{zigzagEQ}
    0\rightarrow A \xrightarrow{i} B \xrightarrow{p} C \rightarrow 0
\end{equation}
where $A,B,C$ is a respective chain complex $C_k(A), C_k(B),C_k(C)$, and so each entry $A,B$ or $C$ can be replaced by a vertical sequence, with vertical arrows being the boundary maps $\partial_k^A,\partial_k^B, \partial_k^C$ from dimension $k$ to $k-1$, with each exact row indexed by $k$. In diagram form, the following diagram with exact rows commute:
\begin{figure}[!htb]
\centering
\begin{tikzcd}\centering
     & \vdots \arrow{d}{\partial^A_{k+1}} & \vdots \arrow{d}{\partial^B_{k+1}} & \vdots \arrow{d}{\partial^C_{k+1}} & \\
    0 \arrow{r} & C_k(A) \arrow{r}{i_k}\arrow{d}{\partial^A_k} & C_{k}(B) \arrow{r}{p_k}\arrow{d}{\partial^B_k} & C_{k}(C)\arrow{r}\arrow{d}{\partial^C_k} & 0 \\
    0 \arrow{r} & C_{k-1}(A)\arrow{r}{i_{k-1}}\arrow{d}{\partial^A_{k-1}} & C_{k-1}(B) \arrow{r}{p_{k-1}}\arrow{d}{\partial^B_{k-1}} & C_{k-1}(C)\arrow{r}\arrow{d}{\partial^C_{k-1}} & 0\\
     & \vdots & \vdots & \vdots & 
\end{tikzcd}
\end{figure}

Then, there is a long exact sequence on homology level:
\begin{center}
\begin{tikzpicture}
\matrix[matrix of nodes,ampersand replacement=\&, column sep=0.5cm, row sep=0.5cm](m)
{
 \& $\cdots$ \& $H_{k+1}(C)$ \\
$H_k(A)$ \& $H_k(B)$ \& $H_k(C)$ \\
$H_{k-1}(A)$ \& $\cdots$ \&  \\
};
\draw[->] (m-1-2) edge (m-1-3)
          (m-1-3) edge[curvedlink] (m-2-1)
          (m-2-1) edge (m-2-2) 
          (m-2-2) edge (m-2-3)
          (m-2-3) edge[curvedlink] (m-3-1)
          (m-3-1) edge (m-3-2);
\filldraw [black] (3.3,-0.4) circle (0.01pt) node[anchor=west] {$d_{k-1,k}$};
\filldraw [black] (3.4,0.9) circle (0.01pt) node[anchor=west] {$d_{k,k+1}$};
\end{tikzpicture}
\end{center}
\end{lem}

Note the only non-trivial maps that arise are $d_{k-1,k}$, and hence the ones we construct through a method called \emph{diagram chasing}.

Take a nontrivial cycle $[c]\in H_k(C)$. This means that there is a nontrivial element $c\in Ker(\partial_k^C)\subseteq C_k(C)$ that is a representative of the equivalence class when quotiented by $Im(\partial^C_{k+1})$. Hence, $\partial_k^C(c)=0$. Exactness of the row \eqref{zigzagEQ} implies that $p$ is surjective, and so there is $b\in C_k(B)$ such that 
\begin{equation}\label{bEQ}
    p_k(b)=c.    
\end{equation}

Moreover, we bring $b$ ``downwards" in the diagram by considering $\partial^B_{k}(b)\in C_{k-1}(B)$. We want to consider $i^{-1}_{k-1}(\partial^B_{k}(b))\in C_{k-1}(A)$, and hence we will have found a candidate for the map $d_{k-1,k}$. For such, note that the diagrams in \eqref{zigzagEQ} commute, and hence 
\begin{equation}\label{comm}
    p_{k-1}(\partial^B_{k}(b))=\partial_k^C(p_k(b)).    
\end{equation}

Due to \eqref{bEQ} and that $c\in Ker(\partial_k^C)$, the equation \eqref{comm} implies that $p_{k-1}(\partial^B_{k}(b))=0$, i.e., $\partial^B_{k}(b)\in Ker(p_{k-1})$. Due to exactness of the row \eqref{zigzagEQ}, $Ker(p_{k-1})=Im(i_{k-1})$, and hence $\partial^B_{k}(b)\in Im(i_{k-1})$. Again, by exactness of the row \eqref{zigzagEQ}, $i$ is injective, and hence there is an unique
\begin{equation}\label{adefZIGZAG}
    a:=i^{-1}_{k-1}(\partial^B_{k}(b))\in C_{k-1}(A).
\end{equation}

Note that $a$ is a nontrivial representative of a cycle $[a]\in H_{k-1}(A)$. Indeed, the commutativity of the diagrams \eqref{zigzagEQ} implies that $i_{k-2}(\partial^A_{k-1}(a))=\partial_{k-1}^B(i_{k-1}(a))$. Due to the definition of $a$ in \eqref{adefZIGZAG}, we obtain that $\partial_{k-1}^B(i_{k-1}(a))=\partial_{k-1}^B(\partial^B_{k}(b))$. Note the latter is zero, since $\partial^B $ is a boundary map, i.e., $\partial^B_{k-1}\circ \partial^B_k=0$. Therefore, we obtain that $i_{k-2}(\partial^A_{k-1}(a))=0$, that is, $\partial^A_{k-1}(a)\in Ker (i_{k-2})$. Lastly, since the row in \eqref{zigzagEQ} is exact, $i$ is injective, which means that $Ker (i_{k-2})=0$, and in turn $\partial^A_{k-1}(a)=0$, meaning that $a$ is a cycle. 

Therefore, we define
\begin{equation}
    d_{k-1,k}([c]):=[a].
\end{equation}

This construction will yield the entries of the connection matrix in Section \ref{sec:matrix}. A dynamical interpretation of this construction is carried in the same section.

\section{The Conley index landscape}\label{sec:index}

For a concise introduction of Conley's theory, see Chapters 22 to 24 in \cite{Smoller83}, and its extension to infinite dimensional systems in \cite{Rybakowski82}. For a short survey, \cite{Mischaikow99}, whereas for a longer exposition, \cite{MischaikowMrozek02}. For the original ideas, follow Conley's monograph \cite{Conley78}. 

Consider the space $\mathcal{T}$ of all topological spaces and the equivalence relation given by $Y\sim Z$ for $Y,Z\in \mathcal{T}$ if, and only if $Y$ is homotopy equivalent to $Z$, that is, there are continuous maps $f:Y\rightarrow Z$ and $g:Z\rightarrow Y$ such that $f\circ g$ and $g\circ f$ are homotopic to $id_Z$ and $id_Y$, respectively. 

An example is when $Y$ is a two dimensional disk, and $Z$ is just a point $y_0$ in the disk. Therefore, we consider $p(y)=y_0$ to be the projection of the disk into the point $y_0$, and $i(y_0)=y_0$ is the inclusion of the point in the disk. Clearly $p\circ i(y_0)=y_0$ is the identity map in $Z$, whereas $L^\tau(y)=(1-\tau)y+\tau y_0$ describes the homotopy between $i\circ p$ and $id_Y$, since $L^0(y)= y$ is the identity map in $Y$, and $L^1(y)= y_0$ is $i\circ p$.

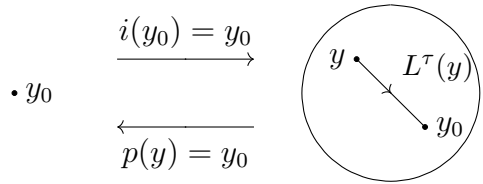
\begin{figure}[ht]\centering
\begin{tikzpicture}[scale=0.9]
    
    \filldraw [black] (-5.5,0) circle (1pt) node[anchor=west]{$y_0$};
    
    \draw[->] (-4,0.5) -- (-2,0.5) node[above] {};
    \filldraw [black] (-3,0.5) circle (0.1pt) node[anchor=south]{$i(y_0)=y_0$};
    \draw[<-] (-4,-0.5) -- (-2,-0.5) node[above] {};
    \filldraw [black] (-3,-0.5) circle (0.1pt) node[anchor=north]{$p(y)=y_0$};
        
    \draw [domain=0:6.28,variable=\t,smooth] plot ({1.3*sin(\t r)},{1.3*cos(\t r)}); 
    \filldraw [black] (0.5,-0.5) circle (1pt) node[anchor=west]{$y_0$};
    \filldraw [black] (-0.5,0.5) circle (1pt) node[anchor=east]{$y$};
    \draw[->] (-0.5,0.5) -- (0,0) ;
    \draw[-] (0,0) -- (0.5,-0.5) ;
    \filldraw [black] (0,0) circle (0.1pt) node[anchor=south west]{\small{$L^\tau(y)$}};    
\end{tikzpicture}
\captionof{figure}{The inclusion of the point $y_0$ into the disk, and the projection of the disk to the point. We also depict the homotopy from the disk to the point $y_0$.}
\end{figure}

Then, the quotient space $\mathcal{T}/\sim$ describes the homotopy equivalent classes $[Y]$ of all topological spaces which have the same homotopy type. Intuitively, $[Y]$ describes all topological spaces which can be continuously deformed into $Y$. In particular, we denote by $[0]$ the trivial homotopy type, that is, the homotopy of a point. Therefore, any contractible space has the homotopy type $[0]$.

The Conley index of a set in phase-space $X$ will be an object in $\mathcal{T}/ \sim$ that captures the topological properties of the local unstable phenomena of a dynamical system.

For any initial data $u_0\in X$, its \emph{orbit} is given by the set its time evolution for all time $\{T(t)u_0\}_{t\in \mathbb{R}}$. Geometrically, an orbit is a curve in phase-space parametrized by time $t$ that passes through the point $u_0$ at initial time $t=0$.

A set $\Sigma\in X$ is \emph{invariant} if $T(t)\Sigma=\Sigma$ for all times $t\in\mathbb{R}$. On one hand, it means that the time evolution of any point in $\Sigma$ still lies in it, and on the other hand, it means that through any point in $\Sigma$ there is a curve globally parametrized for all times that still lies in the set $\Sigma$. Moreover, the structure of invariant sets can change under perturbation, for example when bifurcations occur due to a change of parameter.

Examples of simple invariant sets are equilibria points, heteroclinic or orbits. But invariant sets can also be complicated, such as the Lorenz attractor, Cantor sets coming from Smale horseshoes or chaotic dynamics in general.

We say $\Sigma$ is \emph{isolated} if it has an \emph{isolating neighborhood} $N$, that is, if there is a closed neighborhood $N$ such that $\Sigma$ is contained in the interior of $N$ with $\Sigma$ being the maximal invariant subset of $N$. Note that an isolating neighborhood is robust and remaing isolating with respect to perturbations of the semiflow $T(t)$. In other words, perturbations of the flow does not change the property of being isolating, nevertheless, it can change the dynamics that lie within the invariant set and not in its isolating neighborhood.

From now on, suppose $\Sigma$ is a compact invariant isolated set. We want to understand the dynamics of $\Sigma$ through its associated Conley index. For such, we need to know certain information about the local dynamics nearby $\Sigma$.

The set that captures this local dynamics is $\partial_e N\subset \partial N$ the \emph{exit set of $N$}, that is, the set of points which are not strict ingressing in $N$,
\begin{equation*}
    \partial_e N:=\{u_0\in N \text{ $|$ } T(t)u_0\not\in N \text{ for all sufficiently small $t>0$}\} .
\end{equation*}


\begin{figure}[ht]
\minipage{0.3\textwidth}\centering
\begin{tikzpicture}[scale=0.8]
    \draw [domain=0:6.28,variable=\t,smooth] plot ({1.6*sin(\t r)},{cos(\t r)}) node[anchor=south] {$N$}; 
    
    \draw [line width=0.6mm,domain=1.57:3.14, variable=\t,smooth] plot ({1.6*sin(\t r)},{cos(\t r)}) node[anchor=north]{$\partial_e N$};
    \draw [line width=0.6mm, domain=3.14:4.71,variable=\t,smooth] plot ({1.6*sin(\t r)},{cos(\t r)}) ;   
    
    \draw [domain=-0.5:0.5,variable=\t,smooth] plot ({\t},{(\t)^3}) node[right] {$\Sigma$}; 
\end{tikzpicture}
\endminipage\hfill
\minipage{0.3\textwidth}\centering
\begin{tikzpicture}[scale=0.8]
    \draw [domain=0:6.28,variable=\t,smooth] plot ({1.6*sin(\t r)},{cos(\t r)}) node[anchor= south] {$\partial N\backslash\partial_e N$}; 

    \draw [line width=0.5mm,domain=1.57:3.14, variable=\t,smooth] plot ({1.6*sin(\t r)},{cos(\t r)}) node[anchor=south]{$\partial_e N$};
    \draw [line width=0.5mm, domain=3.14:4.71,variable=\t,smooth] plot ({1.6*sin(\t r)},{cos(\t r)}) ;   
    
    \draw [domain=-0.5:0.5,variable=\t,smooth,->] plot ({\t},{1-(\t)^2}); 
    \filldraw [black] (0,1) circle (0.1pt) node[anchor=north]{$p_1$};

    \draw [domain=-0.5:0.5,variable=\t,smooth,->] plot ({\t},{-1-(\t)^2}); 
    \filldraw [black] (0,-1) circle (1pt) node[anchor=north]{$p_4$};

    \draw[<-] (-1.5,-1) -- (-0.8,-0.5) ;
    \filldraw [black] (-1.1,-0.7) circle (1pt) node[anchor=north]{$p_3$};
    
    \draw[->] (-1.5,1) -- (-0.8,0.5) ;
    \filldraw [black] (-1.1,0.7) circle (0.1pt) node[anchor=north]{$p_2$};
    
    \draw[->] (1,0.55) -- (1.49,0.35) ;
    \draw[->] (1.6,0) -- (2,-0.2) ;
    \draw [domain=0:1.42,variable=\t,smooth,->] plot ({1.6*sin(\t r)},{cos(\t r)}) node[anchor= west] {$\Omega$}; 
    \filldraw [black] (1.6,0) circle (1pt) node[anchor=east]{$p_5$};    
\end{tikzpicture}
\endminipage\hfill
\minipage{0.3\textwidth}\centering
\begin{tikzpicture}[scale=0.8]
    \draw [domain=0:6.28,variable=\t,smooth] plot ({1.6*sin(\t r)},{cos(\t r)}) node[anchor= south] {$\partial N\backslash\partial_e N$}; 

    
    \filldraw [black] (0,-1) circle (1pt) node[anchor=north]{$[\partial_e N]$};

    \draw [domain=-0.5:0.5,variable=\t,smooth] plot ({\t},{(\t)^3}) node[right] {$\Sigma$}; 
\end{tikzpicture}
\endminipage
\captionof{figure}{On the leftmost, we depict an isolated invariant set $\Sigma$, with isolating neighborhood $N$ and exit set $\partial_e  N$. On the middle, the points $p_1,p_2\not\in \partial_e N$, whereas $p_3,p_4\in \partial_e N$. Note that trajectories can drift along $\Omega\subseteq \partial N$, and this drift is not in the exit set, except the exact last time before it leaves, denoted by $p_5$. On the rightmost, we display the quotient of $N$ by $\partial_e N$, that still isolates $\Sigma$.}\label{dinner}
\end{figure}
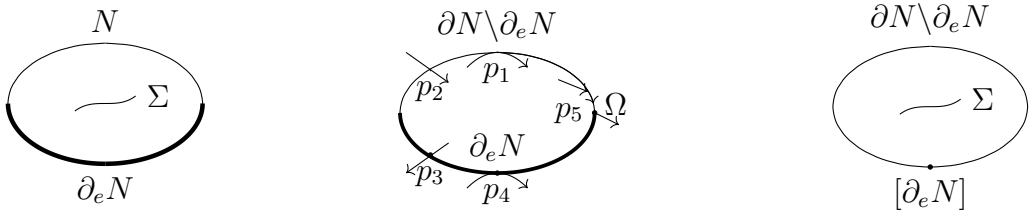

The \emph{Conley index} is defined as
\begin{equation}\label{ConleyIndex}
    h(\Sigma):=[N/ \partial_e N],
\end{equation}
namely the homotopy equivalent class of the quotient space of the isolating neighborhood $N$ relative to its exit set $\partial_e N$. 

Note the quotient space $N/\partial_e N$ has a distinguished point given by $[\partial_e N]$ arising from the quotient topology, which we omit from now on. Also, in case the exit set is empty, we let $h(\Sigma):=[N]\cup \{*\}$, where $*$ denotes a distinguished point outside of $N$ that represents the equivalence class of the empty exit set. The trivial Conley index is given by the homotopy type of a point, and we denote $[0]$.

The usefulness of the index depends on how much information about the structure of the dynamics of the invariant set can be concluded from knowledge of the Conley index. Nevertheless, we mention that the index captures little information on the content of the invariant set itself, but in its surroundings instead: it captures the topology of the unstable region nearby $\Sigma$, and how this set fits within this local dynamics.

\begin{thm}\label{Conleythm}
We now present the three main properties of the Conley index.
\begin{enumerate}
    \item \textbf{\emph{independence of isolating neighborhood:}} if $N$ and $\tilde{N}$ are two isolating neighborhoods of the same invariant set $\Sigma$, then the Conley index is the same: 
    \begin{equation}
        h(\Sigma)=[N/ \partial_e N]=[\tilde{N}/ \partial_e \tilde{N}].
    \end{equation}
    \item \textbf{\emph{additivity:}} if  $\Sigma=\coprod_i \Sigma_i$, then $h(\Sigma)=\coprod_i h(\Sigma_i)$.
    \item \textbf{\emph{homotopy invariance:}} if $N^\tau$ is an isolating neighborhood of $\Sigma^\tau$ for the semiflow $T^\tau:\mathbb{R} \times X\to X$ for all $\tau\in [0,1]$, then the Conley index is constant in $\tau$:
    \begin{equation}
        h(L^0(\Sigma))=h(L^\tau(\Sigma))=h(L^1(\Sigma)).
    \end{equation}
    \item \textbf{\emph{detection of nontrivial dynamics:}} if  $h(\Sigma)\neq [0]$, then $\Sigma\neq\varnothing$.
\end{enumerate}
\end{thm}

We emphasize that $h(\Sigma)$ is an intrinsic property of $\Sigma$, since it does not depend on the neighborhood $N$. The third property is also called \emph{continuation} in the literature. Note that throughout the homotopy, the isolating neighborhood remains the same, whereas the flow is homotoped. This changes both the maxial invariant set within, and the exit set. 
The fourth property is also known in the literature as \emph{Ważeski property}, and it indicates that the index contains information about the dynamics of an invariant set and its relation with any isolating neighborhood around it: knowledge about the dynamics on the boundary of an isolating neighborhood can be enough to determine if the dynamics within is nontrivial. 

The properties from Theorem \ref{Conleythm} resemble the properties of topological degree. Indeed, those are closely related. See McCord \cite{McCord} that proves a generalization of the Poincaré-Hopf index theorem: the sum of the Hopf indices on $\Sigma$ is equal (up to a sign) to the Euler characteristic of the Conley index of the isolated invariant set $\Sigma$.

In general, dealing with the homotopy equivalence classes of
topological spaces is difficult. In turn, the Conley index is hard to compute. Thus to simplify matters we use an algebraic fashion: consider the homology of the pointed topological space $h(\Sigma)$. Mathematically, we consider the \emph{homological Conley index} given by
\begin{equation}
    H(h(\Sigma))=H(N/ \partial_e N)
\end{equation}
where the right hand side $H(N/ \partial_e N)$ is the reduced relative homology for the pair $(N, \partial_e N)$, and $H(N/ \partial_e N)$ is the graded homology group given by $\oplus_{k\in\mathbb{N}_0} H_k(N/ \partial_e N)$.

We compute the Conley index and its homology for five examples.

Firstly, the Conley index of a hyperbolic equilibria $u_+$ with Morse index $n$. Consider $N$ to be a closed ball in $X$ centered at $\Sigma:=\{u_+\}$ with sufficiently small radius such that there is no other equilibria in $N$. Hence $N$ is an isolating neighborhood of $\Sigma$. 

The flow provides a homotopy that contracts along the stable directions to the hyperbolic equilibria $u_+$. Indeed, due to hyperbolicity, the linear flow nearby can be decomposed into stable and unstable linear spaces: $X=E^s \oplus E^u$. Moreover, due to the Hartman-Grobman theorem, the linear flow is locally topologically conjugate to the nonlinear one. Hence, we consider as if we made a change of coordinates, and without loss of generality, we consider the linear flow from now on. 

Furthermore, for any $u_0\in B(u_+)$, we can decompose the linear flow $T(t)u_0$ into stable and unstable parts using the projections $P_u:X\to E^u$ and $P_s:X\to E^s$, yielding $T(t)u_0=P_uT(t)u_0+P_sT(t)u_0$ such that $P_uT(t)u_0\to u_+$ as $t\to -\infty$, and also $P_sT(t)u_0\to u_+$ as $t\to \infty$. 

Therefore, we consider the family of semiflows
\begin{equation}
    T^\tau(t)u_0:=P_uT(t)u_0+P_sT(t+\tan(\tau \pi/2))u_0
\end{equation}
for $\tau\in [0,1]$. This homotopes the original flow $T^0(t)u_0=T(t)u_0$ in the ball $N$ to the flow $T^1(t)u_0=P_uT(t)u_0+u_+$ in the unstable space $N \cap E^u=P_u(N)$. 

Note that the invariant set $\Sigma^\tau$ remains the same throughout the homotopy: it consists only of the equilibria $u_+$. On the other hand, isolating neighborhood $N^\tau$ homotopes $N^0=N$ to an $n$-dimensional ball $N^1=B^n$ in $E^u$. Indeed, $\dim (E^u)=n$, since $u_+$ has Morse index $n$. Therefore $N \cap E^u=P_u(N)$ is a ball $B^n$ of dimension $n$ in $E^u$ with same radius as $N$.

The exit set is now computable, and is given by $\partial_e B^n=\partial B^n$, since after the homotopy there is no more stable direction and the equilibria is hyperbolic. Therefore, everything in its neighborhood leaves after appropriate time. Moreover, the boundary of an $n$ dimensional ball $\partial B^n$ is an $(n-1)$-dimensional sphere $\mathbb{S}^{n-1}$. Therefore, the quotient of a $n$-ball and its boundary is an $n$-sphere. 

\begin{figure}[ht]
\minipage{0.3\textwidth}\centering
    \begin{tikzpicture}[scale=1]
        \draw[->] (-1.5,0) -- (1.5,0) node[right] {$E^s$};
        \draw[<-] (0.49,0) -- (0.495,0);        
        \draw[<-] (-0.49,0) -- (-0.495,0);                

        \draw[->] (0,-1.5) -- (0,1.5) node[right] {$E^u$};
        \draw[->] (0,-0.49) -- (0,-0.495);        
        \draw[->] (0,0.49) -- (0,0.495);                
        
        \draw (0,0) circle (1cm);
        \filldraw [black] (0.6,-0.8) circle (0.1pt) node[anchor=north west] {$N$};

        \filldraw [black] (0,0) circle (1pt) node[anchor=south east]{$u_+$};
    \end{tikzpicture}
\endminipage\hfill
\minipage{0.3\textwidth}\centering%
    \begin{tikzpicture}[scale=1]
        \draw[->] (-2.5,0) -- (-1.5,0) node[anchor=south east] {$\text{homotopy}$};
        \draw[->] (1.5,0) -- (2.5,0) node[anchor=south east] {$\frac{B^n}{\partial_e B^n}$};

        \draw[->] (0,-1.5) -- (0,1.5) node[right] {$E^u$};
        \filldraw [black] (0,0) circle (1pt) node[anchor=east]{$u_+$};

        \draw[->] (0,-0.49) -- (0,-0.495);        
        \draw[->] (0,0.49) -- (0,0.495);                

        \filldraw [black] (0,0) circle (1pt) node[anchor=west]{$B^n$};
        \filldraw [black] (0,1) circle (1pt) node[anchor=west] {$\partial_e B^n$}; 
        \filldraw [black] (0,-1) circle (1pt) ; 
    \end{tikzpicture}
\endminipage\hfill
\minipage{0.3\textwidth}\centering
    \begin{tikzpicture}[scale=1]
        \draw (0,0) circle (1cm);
        \filldraw [black] (0,1) circle (0.1pt) node[anchor=south]{$\mathbb{S}^n$};            
        \filldraw [black] (0,-1) circle (1pt) node[anchor=north]{$[\partial_e N]$};    
    \end{tikzpicture}
\endminipage
\caption{The Conley index of a hyperbolic fixed point $u_+$ with Morse index $n$: homotope the stable directions of the flow, yielding a flow in the unstable subspace $E^u$ of dimension $n$. Then quotient an $n$-dimensional ball with its boundary yields $\mathbb{S}^n$.}
\end{figure}
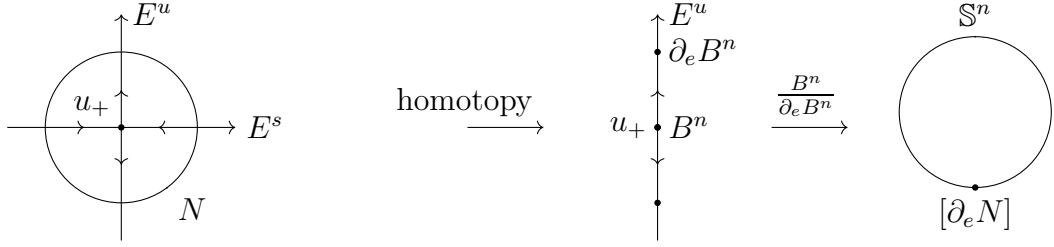

These arguments altogether, using the index is homotopy invariant, yield
\begin{equation}\label{ConleyHypEq}
    h(u_+)=\left[\frac{N}{\partial_e N}\right] = \left[\frac{B^n}{\partial_e B^n}\right] =\left[\frac{B^n}{\mathbb{S}^{n-1}}\right]= \left[\mathbb{S}^n\right].
\end{equation} 


Note that the Morse index $n$ of the equilibrium, which measures how much unstable said equilibrium is, appears as the dimension of the Conley index. therefore, not only the Conley index says how unstable an equilibrium is, but also the local topological behaviour of expanding sets.

We can also compute the homological conley index, yielding
\begin{equation}\label{HypHomindex}
    H_k(h(u_+))=H_k([\mathbb{S}^n])=
    \begin{cases}
        \mathbb{Z}_2  & \text{ if } k=n\\
        0 & \text{ if } k\neq n .
    \end{cases}
\end{equation} 

Secondly, the Conley index of the union of two disjoint invariant sets, for example hyperbolic equilibria $u_-$ and $u_+$ with respective Morse indices $i(u_{-})=m$,$ i(u_{+})=n$, and disjoint isolating neighborhoods $N_-,N_+$. Then, $N:=N_-\cup N_+$ is an isolating neighborhood of $\{ u_-,u_+\}$. 

Since $N_-$ and $N_+$ are disjoint, then $\partial_e (N_-\cup N_+)=\partial_e N_-\cup \partial_e N_+$. By definition of the wedge sum, identifying $\partial_e N_-$ and $\partial_e N_+$, we obtain that 
\begin{align}\label{ConleyTwoEq}
    h(\{ u_-,u_+\})&=\left[\frac{N_-\cup N_+}{\partial_e (N_-\cup N_+)}\right]\\
    &=\left[\frac{N_-}{\partial_e N_-}\vee\frac{N_+}{\partial_e N_+}\right]=h(u_-)\vee h(u_+)
\end{align}
where it was also used that the homotopy type of a wedge sum, is the wedge sum of the homotopy type.

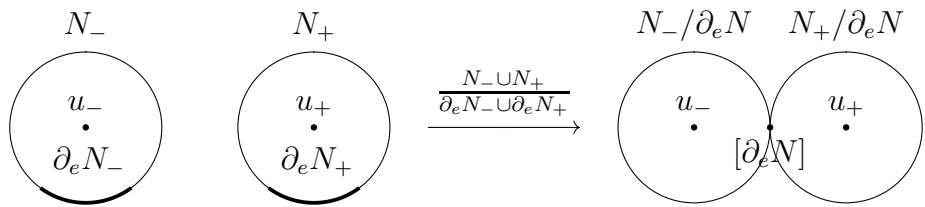
\begin{figure}[ht]\centering
    \begin{tikzpicture}[scale=1]
        \draw (0,0) circle (1cm);
        \filldraw [black] (0,1) circle (0.1pt) node[anchor=south]{$N_-$};            
        \draw [domain=2.5:3.78,line width=0.5mm, variable=\t,smooth] plot ({sin(\t r)},{cos(\t r)}) node[anchor=south west] {$\partial_e N_-$}; 
        \filldraw [black] (0,0) circle (1pt) node[anchor=south]{$u_-$};                
        
        \draw (3,0) circle (1cm);
        \filldraw [black] (3,1) circle (0.1pt) node[anchor=south]{$N_+$};            
        \draw [domain=2.5:3.78,line width=0.5mm, variable=\t,smooth] plot ({3+sin(\t r)},{cos(\t r)}) node[anchor=south west] {$\partial_e N_+$}; 
        \filldraw [black] (3,0) circle (1pt) node[anchor=south]{$u_+$};             
        
        \draw[->] (4.5,0) -- (6.5,0);
        \filldraw [black] (5.5,0) circle (0.1pt) node[anchor=south] {$\frac{N_-\cup N_+}{\partial_eN_- \cup\partial_eN_+}$};

        \draw (8,0) circle (1cm);
        \filldraw [black] (8,1) circle (0.1pt) node[anchor=south]{$N_-/\partial_e N$};            
        \filldraw [black] (8,0) circle (1pt) node[anchor=south]{$u_-$};                
        
        \draw (10,0) circle (1cm);
        \filldraw [black] (10,1) circle (0.1pt) node[anchor=south]{$N_+/\partial_e N$};            
        \filldraw [black] (10,0) circle (1pt) node[anchor=south]{$u_+$};           
        \filldraw [black] (9,0) circle (1pt) node[anchor=north]{$[\partial_e N]$};                     
    \end{tikzpicture}
\caption{The Conley index of two hyperbolic fixed points $u_-,u_+$ with respective Morse indices $m,n$: collapsing the exit sets of its disjoint isolating neighborhoods $N_\pm$ yield a wedge sum of spheres}
\end{figure}

Using the Conley index of a hyperbolic fixed point, as in \eqref{ConleyHypEq}, we obtain that
\begin{equation}
     h(\{ u_-,u_+\})=[\mathbb{S}^m] \vee [\mathbb{S}^n].
\end{equation}

Similarly, we can compute the homological conley index,
\begin{equation}
    H_k(h(\{ u_-,u_+\}))=H_k(\mathbb{S}^m \vee \mathbb{S}^n)=H_k(\mathbb{S}^m )\oplus H_k(\mathbb{S}^n)=
    \begin{cases}
        \mathbb{Z}_2  & \text{ if } k=n,m\\
        0 & \text{ if } k\neq n,m .
    \end{cases}
\end{equation}

Thirdly, we compute the Conley index of a given heteroclinic connection in the plane between two hyperbolic fixed points $u_-$ and $u_+$ of Morse index one. Consider the isolating neighborhood $N$ of the heteroclinic $\Sigma=\{u_\pm\}\cup \{u(t)|u(t)\to_{t\rightarrow_{\pm\infty}} u_\pm\}$ with exit set $\partial_e N$ as depicted in the figure below.

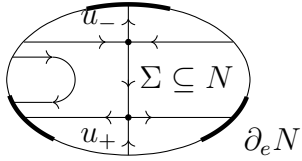
\begin{figure}[ht!]
\centering
\begin{tikzpicture}[scale=1]
    \draw [domain=0:6.28,variable=\t,smooth] plot ({1.6*sin(\t r)},{cos(\t r)}); 
    
    \draw [line width=0.6mm, domain=3.8:4.5,variable=\t,smooth] plot ({1.6*sin(\t r)},{cos(\t r)}) ;   
    \draw [line width=0.6mm,domain=1.8:2.5, variable=\t,smooth] plot ({1.6*sin(\t r)},{cos(\t r)});
    \draw [line width=0.6mm,domain=-0.35:0.35, variable=\t,smooth] plot ({1.6*sin(\t r)},{cos(\t r)});
    
    \filldraw [black] (0,-0.5) circle (1pt) node[anchor=north east]{$u_{+}$};
    \draw (-1.375,0.5) -- (1.375,0.5);    
    \draw[<-] (0.3,0.5) -- (0.31,0.5);   
    \draw[<-] (-0.3,0.5) -- (-0.31,0.5);   
    \draw[->] (0,0.8) -- (0,0.81);   
    
    \filldraw [black] (0,0.5) circle (1pt) node[anchor=south east]{$u_{-}$};    
    \draw [domain=-1:1,variable=\t,smooth] plot ({0},{(-\t)});
    \draw[<-] (0,-0.1) -- (0,0)    
    node[right] {$\Sigma\subseteq N$};   
    \draw (-1.375,-0.5) -- (1.375,-0.5) node[anchor=north west]{$\partial_e N$};
    \draw[->] (-0.3,-0.5) -- (-0.31,-0.5);   
    \draw[->] (0.3,-0.5) -- (0.31,-0.5);   
    \draw[<-] (0,-0.8) -- (0,-0.81);   

    \draw[->] (-1.525,0.3) -- (-1,0.3);    
    \draw (-1.525,-0.3) -- (-1,-0.3);  
    \draw [->,domain=0:3.14, variable=\t,smooth] plot ({-1+0.3*sin(\t r)},{0.3*cos(\t r)});    
\end{tikzpicture}
\end{figure}

We identify the exit set consecutively, yielding $[\partial_e N]$ as a point in the quotient topology:

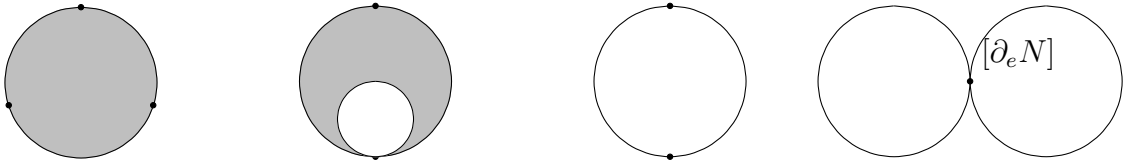
\begin{figure}[H]
\minipage{0.25\textwidth}\centering
\begin{tikzpicture}[scale=1]
    \filldraw [lightgray, domain=0:6.28,variable=\t,smooth] plot ({sin(\t r)},{cos(\t r)}); 
    \draw [domain=0:6.28,variable=\t,smooth] plot ({sin(\t r)},{cos(\t r)}); 

    \filldraw [black] (0,1) circle (1pt);
    \filldraw [black] (-0.95,-0.3) circle (1pt);
    \filldraw [black] (0.95,-0.3) circle (1pt);
\end{tikzpicture}
\endminipage\hfill
\minipage{0.25\textwidth}\centering
\begin{tikzpicture}[scale=1]
    \filldraw [lightgray, domain=0:6.28,variable=\t,smooth] plot ({sin(\t r)},{cos(\t r)}); 
    \draw [domain=0:6.28,variable=\t,smooth] plot ({sin(\t r)},{cos(\t r)}); 
    
    \filldraw [black] (0,1) circle (1pt);
    \filldraw [black] (0,-1) circle (1pt);

    \filldraw [white, domain=0:6.28,variable=\t,smooth] plot ({0.5*sin(\t r)},{-0.5+0.5*cos(\t r)}); 
    \draw [domain=0:6.28,variable=\t,smooth] plot ({0.5*sin(\t r)},{-0.5+0.5*cos(\t r)});
\end{tikzpicture}
\endminipage\hfill
\minipage{0.25\textwidth}\centering
\begin{tikzpicture}[scale=1]
    \draw [domain=0:6.28,variable=\t,smooth] plot ({sin(\t r)},{cos(\t r)}); 
    
    \filldraw [black] (0,1) circle (1pt);
    \filldraw [black] (0,-1) circle (1pt);
\end{tikzpicture}
\endminipage\hfill
\minipage{0.25\textwidth}\centering
\begin{tikzpicture}[scale=1]
    \draw [domain=0:6.28,variable=\t,smooth] plot ({sin(\t r)},{cos(\t r)}); 
    \draw [domain=0:6.28,variable=\t,smooth] plot ({2+sin(\t r)},{cos(\t r)}); 
    
    \filldraw [black] (1,0) circle (1pt) node[anchor=south west]{$[\partial_e N]$};
\end{tikzpicture}
\endminipage
\captionof{figure}{Consider the isolating neighborhood $N$ of $\Sigma=\{u_i\}\cup \mathcal{C}_{i-1,i}\cup u_{i-1}$ with exit set $\partial_e N_I$ such that $N\subseteq N_I$ is an isolating neighborhood of $M_{i-1}$ with $\partial_e N=\partial_e N_I$. On the middle figure, we consider the quotient of the exit set set, which does not affect $\Sigma$ since it is an isolating neighborhood. On the rightmost drawing, we see the quotient of $N_I$ by $N$, which is an isolating neighborhood of $M_i$ with exit set $[N]$. Note that in the latter case, the exit set consists only of orbits from the unstable manifold of $M_i$.}
\end{figure}

Therefore, the Conley index of such heteroclinic connection between two equilibria points of saddle type is the wedge sum of two spheres - the same as the Conley index of two isolated points. Therefore, the Conley index loses dynamical information of the isolated invariant set. In other words, one can not infer the dynamical structure of the isolated invariant set, given its Conley index. The index simply carries information about the instability of the isolated invariant set in its surroundings, not about itself.

Now, we compute the Conley index of a normally hyperbolic periodic orbit given by $\Sigma:=\{ u_0\in X \text{ $|$ } T(t)u_0=T(t+p)u_0\}$ such that its associated Poincaré map $P$ has derivative $DP$ with $n$ positive eigenvalues bigger than 1. Consider the tubular neighborhood $N_\epsilon:=\{v\in X  : ||v-u_0|| \leq \epsilon \text{ for } u_0\in \Sigma\}$ which is isolating, for sufficiently small $\epsilon>0$. 

Note that due to normal hyperbolicity, its normal directions (which correspond to the Poincaré section) can be split into stable and unstable directions, whereas its tangent direction traces nearby points of the periodic orbit yielding the Poincaré map. Mathematically, $X=E^s \oplus E^{cu}$ where $E^s$ is the stable subspace, and $E^{cu}$ correspond to the center-unstable subspace, i.e., the unstable space and the tangent space of $\Sigma$. Similarly as the case of hyperbolic fixed point, we homotope the stable directions of the flow. Therefore, after homotopy, we can consider $N^{cu}_\epsilon$ to be the tubular neighborhood of $\Sigma$ in $E^{cu}$. Note that $N^{cu}$ isomorphic to $\mathbb{S}^1\times B^n_\epsilon$ where $\mathbb{S}^1$ parametrizes the angle of the periodic orbit, and $B^n_\epsilon$ is the ball of dimension $n$ and radius $\epsilon$ in the strict unstable subspace $E^u$.

The exit set is given by $\partial_e N^{cu}_\epsilon=\partial_e (\mathbb{S}^1\times B^n_\epsilon)=\mathbb{S}^1\times \partial B^n_\epsilon=\mathbb{S}^1\times \mathbb{S}^{n-1}$, since after the homotopy there is no more stable direction of the normally hyperbolic periodic orbit, and for each point in the periodic orbit $\Sigma$ component, the exit set consists of the whole boundary $\partial B^n$. Moreover, the boundary of an $n$ dimensional ball $\partial B^n$ is an $(n-1)$-dimensional sphere $\mathbb{S}^{n-1}$. 

\begin{figure}[ht]
\minipage{0.3\textwidth}\centering
    \begin{tikzpicture}[scale=1]
        
        \draw (0,0) circle (0.6cm);
        \draw[line width=0.5mm] (0,0) circle (1cm);
        \draw[line width=0.5mm] (0,0) circle (0.2cm);

        \filldraw [black] (0,1) circle (0.1pt) node[anchor=south] {$\partial_e N$};
        \draw[->] (0.6,-0.01) -- (0.6,0) node[anchor=south] {$\Sigma$};
    \end{tikzpicture}
\endminipage\hfill
\minipage{0.3\textwidth}\centering%
    \begin{tikzpicture}[scale=1]

        \draw (-1,0) arc (180:360:1cm and 0.5cm);
        \draw[dashed] (-1,0) arc (180:0:1cm and 0.5cm);
        \draw[->] (-0.1,-0.5) -- (0,-0.5) node[anchor=south] {$\Sigma$};
        
        \draw (0,0) circle (1cm);
    
        \filldraw [black] (0,1) circle (1.5pt) node[anchor=south] {$[\partial_e N]$}; 
        \filldraw [black] (0,-1) circle (1.5pt) ; 
    \end{tikzpicture}
\endminipage\hfill
\minipage{0.3\textwidth}\centering
    \begin{tikzpicture}[scale=1]
        \draw (0,0) circle (1cm);
        \filldraw [black] (0,1) circle (0.1pt) node[anchor=south]{$PT^2$};
        \draw (0.5,0) circle (0.5cm);        
        \filldraw [black] (1,0) circle (1.5pt) node[anchor=west]{$[\partial_e N^{cu}_\epsilon]$};
        
        \draw (-1,0) arc (180:360:0.5cm and 0.25cm);
        \draw[dashed] (-1,0) arc (180:0:0.5cm and 0.25cm);
        \draw[->] (-0.54,-0.25) -- (-0.5,-0.25) node[anchor=north] {$\Sigma$};        
    \end{tikzpicture}
\endminipage
\caption{The index of a periodic orbit with Morse index $1$: after a homotopy of the flow along the stable directions, we obtain an annulus $N\subseteq E^{cu}$. Identifying the outer ring (obtaining the north pole), and inner ring (resulting in the south pole), we obtain a sphere such that the equator is the periodic orbit. Then, we identify the north and south pole, yielding the pinched torus.}\label{pinched}
\end{figure}
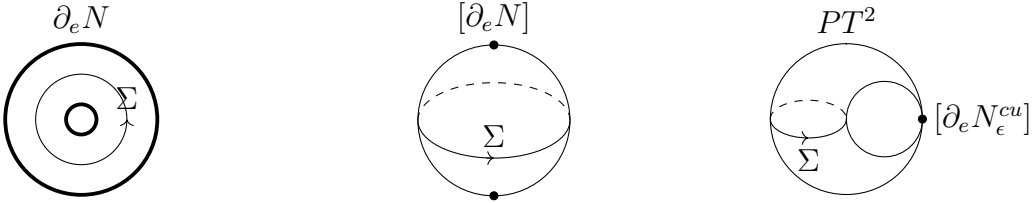

These arguments altogether, using the index is homotopy invariant, yield
\begin{equation}\label{ConleyPeriodic}
    h(\Sigma)=\left[\frac{N_\epsilon}{\partial_e N_\epsilon}\right] =\left[\frac{N_\epsilon^{cu}}{\partial_e N_\epsilon^{cu}}\right] =\left[\frac{\mathbb{S}^1\times B^n_\epsilon}{ \mathbb{S}^1\times \mathbb{S}^{n-1}}\right]=[PT^{n+1}]
\end{equation} 
where $PT^{n+1}$ denotes the skewed $(n+1)$-dimensional pinched torus. We explain how to obtain the pinched torus from the quotient $\mathbb{S}^1\times B^n_\epsilon/ \mathbb{S}^1\times \mathbb{S}^{n-1}$ in Figure \ref{pinched}, for the case $n=1$. For $n>1$, 
one needs to identify firstly each $B^n_\epsilon/\mathbb{S}^{n-1}=\mathbb{S}^{n}$, and then identify one point from each of those $n$-dimensional sphere.

Therefore its homological Conley index is
\begin{equation}\label{PTHomindex}
    H_k(PT^{n+1})=
    \begin{cases}
        \mathbb{Z}_2  & \text{ if } k=n,n+1 \\
        0 & \text{ if } k\neq n,n+1 .
    \end{cases}    
\end{equation}

In general, given an isolated invariant set, we might be able to compute its Conley index depending on how complicated such set and its isolating neighborhood are.  

\begin{figure}[ht]
\minipage{0.3\textwidth}\centering
    \begin{tikzpicture}[scale=0.9]
        
        \draw (0,0) circle (0.6cm);
        \draw[line width=0.5mm] (0,0) circle (1cm);
        \draw[line width=0.5mm] (0,0) circle (0.2cm);

        \filldraw [black] (0,1) circle (0.1pt) node[anchor=south] {$\partial_e N$};
        \draw[->] (0.6,-0.01) -- (0.6,0);
    \end{tikzpicture}
\endminipage\hfill
\minipage{0.3\textwidth}\centering%
    \begin{tikzpicture}[scale=0.9]
        \draw (0,0) circle (0.6cm);
        \draw[line width=0.5mm] (0,0) circle (1cm);
        \draw[line width=0.5mm] (0,0) circle (0.2cm);

        \filldraw [black] (0,-0.6) circle (1pt) node[anchor=south] {$u_*$};
        \filldraw [black] (0,1) circle (0.1pt) node[anchor=south] {$\partial_e N$};
        \draw[->] (0.6,-0.01) -- (0.6,0);
    \end{tikzpicture}
\endminipage\hfill
\minipage{0.3\textwidth}\centering
    \begin{tikzpicture}[scale=0.9]
        \draw (0,0) circle (0.6cm);
        \draw[line width=0.5mm] (0,0) circle (1cm);
        \draw[line width=0.5mm] (0,0) circle (0.2cm);

        \filldraw [black] (0,-0.6) circle (1pt) node[anchor=south] {$u_*$};
        \filldraw [black] (0,0.6) circle (1pt) node[anchor=north] {$u_{**}$};        
        \filldraw [black] (0,1) circle (0.1pt) node[anchor=south] {$\partial_e N$};
        \draw[->] (0.6,-0.01) -- (0.6,0);         \draw[->] (-0.6,0.01) -- (-0.6,0);       
    \end{tikzpicture}
\endminipage
\caption{The Conley index of a normally hyperbolic periodic orbit, homoclinic orbit, and heteroclinic cycle are the same if they have the same isolating neighborhood.}\label{portugal}
\end{figure}
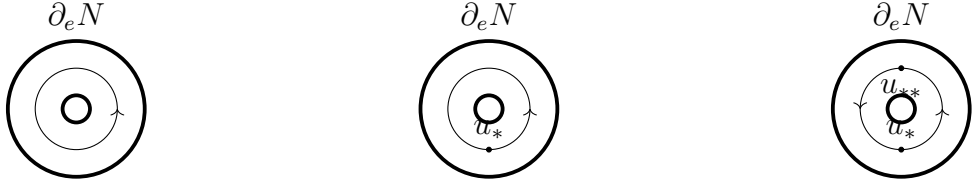

Conversely, given a Conley index, it is not possible to conclude the structure of the isolated invariant set. This occurs due to the homotopy invariance, that opens doors for the occurrence of bifurcations within the invariant set that pass unnoticed in the nearby unstable dynamics. For instance, instead of a periodic orbit in the above example, we could have a homoclinic orbit, or a heteroclinic cycle, as in Figure \ref{portugal}. In all these cases, those sets have the same isolating neighborhood, yielding the same index. Therefore, how to identify the structure of an invariant set, given its index, remains an open question.

In some cases the Conley index with extra assumptions might yield the flow structure in the isolated invariant set. It was proven by McCord, Mrozek and Mischaikow \cite{McCordMrozekMischaikow95} that if the Conley index is a pinched torus as above, and the isolating neighborhood has a Poincaré section, then there is a periodic orbit in the isolating neighborhood. 

As a last computable example, we suppose that the semiflow $T(t)u_0$ has a compact global attractor $\mathcal{A}$ in a Banach space $X$, and we compute its Conley index. Two different proofs of this fact are in Corollary 3.2 of \cite{McCordMischaikow96}, or Theorem 6.2 in \cite{HattoriMischaikow91}. 

Note there is a sufficiently large ball $B$ in $X$ with empty exit set. Otherwise, if any ball of sufficiently big radius have nonempty exit set, then we can construct a solution which is not bounded, a contradiction to dissipativity. Therefore, we can compute the Conley index from its definition, yielding
that the index of $\mathcal{A}$ has the same homotopy type of the ball, which is the homotopy type of a point:
\begin{equation}\label{ConleyAtt}
    H_k(h(\mathcal{A})))=
    \begin{cases}
        \mathbb{Z}_2  & \text{ if } k=0\\
        0 & \text{ if } k\neq 0.
    \end{cases}
\end{equation}

Notice the Conley index of the global attractor is the same as a globally stable fixed point. This means that the evolution of any neighborhood of the attractor can not leave and stay away from the attractor too much, but it has to converge to the attractor In other words, there is a suitable neighborhood such that its exit set is empty. 

Similarly to the the stable fixed point, note that the cellular homology of the Conley index has $\mathbb{Z}_2\oplus \mathbb{Z}_2$ in dimension zero, one component regarding the fixed point, and the other regarding the added distinguished point in order to represent the empty exit set. In reduced homology, we quotient this by $\mathbb{Z}_2$ in order to obtain \eqref{ConleyAtt}.

We emphasize once again: the knowledge of the Conley index of the attractor does not impose any restriction on the topology of global attractors, but on its surroundings only: the global attractor is an asymptotically global stable set. Nevertheless, since the exit set is empty, then we know that $[\mathcal{A}]$ has the homotopy type of a point. 

Lastly, the Conley index is not always computable. There are sets which do not admit an isolating neighborhood. Take as an example a center fixed point, that is, in any local neighborhood, there are periodic orbits winding around the fixed point. This center fixed point can not be isolated, that is, there is no neighborhood that has a maximal invariant set such that its closure is contained in the interior of the neighborhood.

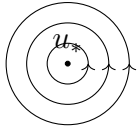
\begin{figure}[ht]
\centering
    \begin{tikzpicture}[scale=0.9]
        \draw (0,0) circle (0.6cm);
        \draw (0,0) circle (0.9cm);
        \draw (0,0) circle (0.3cm);
        \filldraw [black] (0,0) circle (1pt) node[anchor=south] {$u_*$};
        
        \draw[->] (0.6,-0.01) -- (0.6,0);         
                \draw[->] (0.9,-0.01) -- (0.9,0);         
                        \draw[->] (0.3,-0.01) -- (0.3,0); 
    \end{tikzpicture}
\caption{A center equilibrium $u_*$ does not admit an isolating neighborhood.}\label{portugalHUNGRY}
\end{figure}

\section{A tool called connection matrix}\label{sec:matrix}

The general framework for the Conley index and the connection matrix rely on the decomposition of an invariant set $\Sigma$ into smaller invariant sets. That is, subsets of $\Sigma$ that have the dynamical characteristic of attracting or repelling. Any other point in $\Sigma$ is a heteroclinic connection between attractive and repulsive sets. In particular, those attractive or repulsive sets can sometimes be indexed, and carry an order. This is called a Morse decomposition. See \cite{Milnor63} for mathematical definition. Moreover, the connections between Morse sets can be computed by relating the Conley indices of Morse sets to the Conley index of the total isolated invariant set they are contained. In particular, the connection problem is reduced to a problem in linear algebra. More precisely we will show that given an isolated invariant set with certain attracting and repelling subsets, the relationship between their Conley indices can be expressed in terms of a matrix called the connection matrix.

Since the semiflow $T(t)$ is gradient, there is a Lyapunov function $L:X\to \mathbb{R}$ which is decreasing along orbits. We suppose $L$ does not decrease for finitely many energy levels $l_i\in \mathbb{R}$ with $i=0,...,N$. We consider the sets $M_i$ to be the maximal compact invariant sets in $L^{-1}(l_i)$. Therefore, the semiflow $T(t)$ is of steepest descent, which mimics the decay in the induced decay semiflow in $\mathbb{R}$ through $L$, except at $\mathcal{R}=\cup_{i=0}^N M_i$.

Hence, $\mathcal{M}:=\{M_i\}_{i=0}^N$ is a \emph{Morse decomposition} of the global attractor $\mathcal{A}$: the \emph{Morse sets} $M_i$ are mutually disjoint and invariant subsets of $\mathcal{A}$ such that for any $u_0\in\mathcal{A}$, either $u_0\in \mathcal{M}$, or there exists a global solution $T(t)u_0$ such that 
\begin{equation}\label{het}
    M_i\xleftarrow{t\to-\infty} T(t)u_0\xrightarrow{t\to\infty} M_j \text{ for } i>j.
\end{equation}

The Morse set $M_0$ is called the \emph{stable Morse set}, since there is no trajectory leaving it within the attractor, due to \eqref{het}. We will consider the \emph{homology of Morse sets} given by the homology of the Conley index of such sets, namely $H(h(M_i))$.

If $u_0$ is an initial condition with global solutions satisfying \eqref{het}, we call it a \emph{heteroclinic solution} that connects $M_i$ to $M_j$. We denote such set of connections by $\mathcal{C}_{i,j}=\mathcal{C}(M_i,M_j)$. The goal of this section is to construct a matrix $\Delta$ with entries $\Delta_{ji}$ that will tells us if $C_{i,j}\neq 0$.


The \emph{connection matrix} is a linear map given by
\begin{equation}\label{Delta}
    \Delta: H(h(\mathcal{M}))\to H(h(\mathcal{M}))
\end{equation}
where $H(h(\mathcal{M}))=\oplus_{i=0}^N H(h(M_i))$, since the sets $M_i$ are disjoint. 

Therefore, $\Delta$ can be represented as the following $(N+1)\times (N+1)$ matrix form
\begin{equation}\label{DeltaMatrix}
  \Delta=
  \begin{blockarray}{*{3}{c} l}
    \begin{block}{*{3}{>{$\footnotesize}c<{$}} l}
      $H(h(M_0))$ & $\cdots$ & $H(h(M_{N}))$ \\
    \end{block}
    \begin{block}{[*{3}{c}]>{$\footnotesize}l<{$}}
      \Delta_{00} & \cdots & \Delta_{0N} \bigstrut[t] & $H(h(M_0))$ \\
      \vdots & \ddots & \vdots & $\vdots$ \\
      \Delta_{N0} & \cdots & \Delta_{NN} & $H(h(M_N))$. \\
    \end{block}
  \end{blockarray}
\end{equation}
Each entry of \eqref{DeltaMatrix} is given by 
\begin{equation}\label{Deltaij}
    \Delta_{ij}: H(h(M_j))\to H(h(M_i))
\end{equation}
where $H(h(M_i))=\oplus_{i^{'}=1}^{C_i} H(h(M_{ii^{'}}))$, since each Morse set $M_i$ has $C_i$ connected compontents denoted by $M_i=\cup_{i^{'}=1}^{C_i} M_{i,i^{'}}$.

Hence, the entries $\Delta_{ij}$ of the matrix \eqref{DeltaMatrix} are still 
$C_j\times C_i$ matrices given by
\begin{equation}\label{DeltaijMatrix}
  \Delta_{ij}=
  \begin{blockarray}{*{3}{c} l}
    \begin{block}{*{3}{>{$\footnotesize}c<{$}} l}
      $H(h(M_{j,1}))$ & $\cdots$ & $H(h(M_{j,C_j}))$ \\
    \end{block}
    \begin{block}{[*{3}{c}]>{$\footnotesize}l<{$}}
      \delta_{11} & \cdots & \delta_{1C_j} \bigstrut[t] & $H(h(M_{i,1}))$ \\
      \vdots & \ddots & \vdots & $\vdots$ \\
      \delta_{C_i1} & \cdots & \delta_{C_iC_j} & $H(h(M_{i,C_i}))$. \\
    \end{block}
  \end{blockarray}
\end{equation}

Furthermore, each entry of \eqref{DeltaijMatrix} is a linear map given by
\begin{equation}\label{deltaij}
    \delta_{i'j'}: H(h(M_{jj'}))\to H(h(M_{ii'}))    
\end{equation}
where the graded homology can be decomposed in each homology level as $H(h(M_{ii'})) = \oplus_{k=0}^{N_{i'}} H_k(h(M_{ii'}))$ such that $N_{i'}$ denote the highest nontrivial Homology level for $h(M_{ii'})$.

Finally, the entries $\delta_{i'j'}$ of $\Delta_{ij}$ are $(N_j+1)\times (N_i+1)$ matrices as well:
\begin{equation}\label{deltaijMatrix}
  \delta_{i'j'}=
  \begin{blockarray}{*{3}{c} l}
    \begin{block}{*{3}{>{$\footnotesize}c<{$}} l}
      $H_0(h(M_{jj'}))$ & $\cdots$ & $H_{N_{j'}}(h(M_{jj'}))$ \\
    \end{block}
    \begin{block}{[*{3}{c}]>{$\footnotesize}l<{$}}
      d_{00} & \cdots & d_{0N_{j'}} \bigstrut[t] & $H_0(h(M_{ii'}))$ \\
      \vdots & \ddots & \vdots & $\vdots$ \\
      d_{N_{i'}0} & \cdots & d_{N_{i'}N_{j'}} & $H_{N_{i'}}(h(M_{ii'}))$ \\
    \end{block}
  \end{blockarray}
\end{equation}
where $d_{mn}:H_n(h(M_{jj'}))\to H_m(h(M_{ii'}))$ is a linear map between vector spaces.






Consider the interval $I:=\{i^*,i^*+1,...,j^*-1,j^*\}$ that consists of all indices from $i^*$ up to $j^*$. Definte the sets $M_I:=\cup_{i\in I} M_i \cup_{i,j\in I}\mathcal{C}_{i,j}$. The homology of the Conley index of the sets $M_I$ will be called \emph{homology index braid}, denoted by $H(h(M_I))$. 

Moreover, we denote by $\Delta_I$ to be a 
sub-matrix of $\Delta$ having entries $\Delta_{i,j}$ given by \eqref{DeltaijMatrix} with $i,j\in I$.

The main result of existence of the matrix $\Delta$ that encodes the entanglement of Morse set is due to Franzosa \cite{Franzosa89}. Below are its main properties. 

\begin{thm}\label{ThmDelta}
There exists a matrix $\Delta$ as \eqref{Delta} with the following properties:
\begin{enumerate}
    \item \textbf{\emph{strict upper triangular:}} if $i\geq j$, then $\Delta_{i,j}=0$.
    \item \textbf{\emph{boundary operator:}} $\Delta^2=0$.
    \item \textbf{\emph{degree -1 map:}} $\Delta (H_k(h(M_{jj'})))\subseteq H_{k-1}(h(M_{ii'}))$ for some $i$. 
    \item \textbf{\emph{homology map:}} $H(h(M_I))\cong Ker \Delta_I / Im \Delta_I$ for any interval $I$.
\end{enumerate}
\end{thm}


The first item is a reflection that the order of the Morse decomposition is preserved: the only possible connections occur from bigger $j$ to lower $i$. 

The second item, shows that indeed $\Delta$ is a boundary operator in the topological sense, as in Section \ref{sec:alg}. This implies that there entries satisfy a linear system, and hence the entries of this matrix are entangled and dependent on one another. 

The third item, shows that the boundary map actually maps a $k^{th}$-level homology to $(k-1)^{th}$-level. In other words, $d_{mn}=0$ if $m\neq n-1$, i.e., it is possibly nonzero only when it is a map from $H_n(h(M_{jj'})))$ to $H_{n-1}(h(M_{ii'}))$. That is, the only possibly nonzero entries of \eqref{deltaijMatrix} are the upper-diagonal ones, given by $d_{n-1,n}$. This limits the connection matrix tool to detect only heteroclinics when the Morse index of hyperbolic equilibria differ by one. Nevertheless, there is an effort into understanding spectral sequences of a Morse filtration of the attractor, yielding such connections of equilibria with Morse index differing by more than one, as in Theorem 6.1 of \cite{CorneaRezendeSilveira10}. 

The fourth item says that any sub-matrix $\Delta_I$ of $\Delta$ is a homology map compatible with the homology index braid, or conversely: the homology index braid can be computed through sub-matrices of $\Delta$. 

The main consequence is that we know which Morse set connects to which other through a heteroclinic orbit, whenever we can compute the connection matrix.

\begin{thm}\emph{\textbf{Heteroclinic detection.}}\label{thmdetection}
If $\Delta_{i,i+1}\neq 0$, then $\mathcal{C}_{i+1,i}\neq\varnothing$. More precisely, if the entry $\delta_{i^{'}j^{'}}$ of $\Delta_{i,i+1}$ is nonzero, then there is a heteroclinic from the $j^{'}$ connected component of $M_{i+1}$ to the $i^{'}$ connected component of $M_{i}$.
\end{thm}


We now construct the entries of the connection matrix which are going to be useful due to Theorem \ref{thmdetection} in order to detect heteroclinic orbits $\mathcal{C}_{i+1,i}$ between adjacent Morse sets $M_{i+1}$ and $M_i$. For such, consider $I=\{i,i+1\}$, and hence $M_I:=M_{i}\cup M_{i+1} \cup\mathcal{C}_{i+1,i}$. Therefore, the relevant entry in \eqref{DeltaMatrix} for adjacent connections is $\Delta_{i,i+1}$. For the sake of simplicity, we suppose that both $M_{i+1}$ and $M_i$ have a single connected component, and hence the matrix $\Delta_{i,i+1}$ consists of a single entry, which we denote by $\delta_{i,i+1}$.

Consider an isolating neighborhood $N_I$ of $M_I$ with exit set $\partial_e N_I$, and $N\subseteq N_I$ that isolates $M_{i}$ such that $\partial_e N=\partial_e N_I$. See \cite{Kurland} for the existence of such sets. 
In our particular case, let $J_I:=[l_i-\epsilon,l_{i+1}+\epsilon]$ with $\epsilon>0$, and consider its preimage through the Lyapunov function, $N_I:=L^{-1}(J_I)$. Similarly, let $J:=[l_i-\epsilon,l_{i+1}-\epsilon]$ in order to obtain $N:=L^{-1}(J)$. Note that $\partial_e N=\partial_e N_I$, since a solution that exits $N_I$ must decrease its energy below the level $l_{i+1}$, and hence lies in $N$ for sufficiently small $\epsilon>0$.

The inclusion $N\hookrightarrow N_I$ and projection $N_I\rightarrow N_I/N$ induce the quotient maps of inclusion $inc:N/\partial_e N\to N_I/\partial_e N$ and projection $proj:N_I/\partial_e N\to N_I/N$ since $\partial_e N\subseteq N$. Moreover, we obtain the following short exact sequence
\begin{equation}
    0\rightarrow \frac{N}{\partial_e N}\xrightarrow{inc}\frac{N_I}{\partial_e N_I}\xrightarrow{proj}\frac{N_I}{N} \rightarrow 0
\end{equation}
where $proj \circ inc=0$, that is, $Im(inc)=Ker(proj)$.

\begin{figure}[ht]
\minipage{0.33\textwidth}\centering
\begin{tikzpicture}[scale=1]
    \draw [domain=0:6.28,variable=\t,smooth] plot ({sin(\t r)},{1.6*cos(\t r)}) node[anchor=south] {$N_I$}; 
    \draw [fill=lightgray,domain=1.45:4.85,variable=\t,smooth] plot ({sin(\t r)},{1.6*cos(\t r)}) ;   
    \draw [line width=0.6mm, domain=3.14:4.4,variable=\t,smooth] plot ({sin(\t r)},{1.6*cos(\t r)}) ;   
    \draw [line width=0.6mm,domain=1.88:3.14, variable=\t,smooth] plot ({sin(\t r)},{1.6*cos(\t r)}) node[anchor=north]{$\partial_e N_I=\partial_e N$};
    
    \filldraw [black] (0,-1) circle (1pt) node[anchor=west]{\tiny{$M_{i}$}};
    \filldraw [black] (0,1) circle (1pt) node[anchor=west] {\tiny{$M_{i+1}$}};    
    \draw [domain=-1:1,variable=\t,smooth] plot ({0},{(-\t)});
    \draw[<-] (0,-0.1) -- (0,0)    node[right] {\tiny{$M_I$}};   

    \draw (-1,0.2) -- (1,0.2);
    \filldraw [black] (-1.2,0.2) circle (0.01pt) node[anchor=north]{$N$};       
\end{tikzpicture}
\endminipage\hfill
\minipage{0.33\textwidth}\centering
\begin{tikzpicture}[scale=1]
    \draw [domain=0:6.28,variable=\t,smooth] plot ({sin(\t r)},{1.6*cos(\t r)}) node[anchor=south] {$N_I/\partial_e N_I$}; 
    \draw [fill=lightgray,domain=1.45:4.85,variable=\t,smooth] plot ({sin(\t r)},{1.6*cos(\t r)}) ;   
    \filldraw [black] (0,-1.6) circle (1.5pt) node[anchor=north]{$[\partial_e N_I]=[\partial_e N]$};
    
    \filldraw [black] (0,-1) circle (1pt) node[anchor=west]{\tiny{$M_{i}$}};
    \filldraw [black] (0,1) circle (1pt) node[anchor=west] {\tiny{$M_{i+1}$}};    
    \draw [domain=-1:1,variable=\t,smooth] plot ({0},{(-\t)});
    \draw[<-] (0,-0.1) -- (0,0)    node[right] {\tiny{$M_I$}};   

    \draw (-1,0.2) -- (1,0.2);
    \filldraw [black] (-1.4,0.2) circle (0.01pt) node[anchor=north] {$\frac{N}{\partial_e N}$};      
\end{tikzpicture}
\endminipage\hfill
\minipage{0.33\textwidth}\centering
\begin{tikzpicture}[scale=1]
    \draw [domain=0:6.28,variable=\t,smooth] plot ({sin(\t r)},{1.6*cos(\t r)}) node[anchor=south] {$N_I/N$}; 
    \filldraw [black] (0,-1.6) circle (1.5pt) node[anchor=north]{$[N]$};

    \filldraw [black] (0,1) circle (1pt) node[anchor=west] {\tiny{$M_{i+1}$}};    
    \draw [domain=-1:1.6,variable=\t,smooth] plot ({0},{(-\t)});
    \draw[<-] (0,-0.1) -- (0,0)    node[right] {\tiny{$\frac{M_I}{N}$}};   
\end{tikzpicture}
\endminipage
\captionof{figure}{On the leftmost, we see the isolating neighborhood $N_I$ of $M_I$ with exit set $\partial_e N_I$ such that $N\subseteq N_I$ isolates $M_{i}$ with $\partial_e N=\partial_e N_I$. On the middle, we consider the quotient of the exit set set, which does not affect $\Sigma$ since it is an isolating neighborhood. On the rightmost, we see the quotient of $N_I/N$, which is an isolating neighborhood of $M_{i+1}$ with exit set $[N]$. Note that in the latter case, the exit set consists only of orbits from the unstable set of $M_{i+1}$.}
\end{figure}

The maps $inc,proj$ induces maps in chain complex levels. Therefore, due to the zig-zag lemma \ref{thmzigzag}, there is a long exact sequence on homology level:
\begin{center}
\begin{tikzpicture}
\matrix[matrix of nodes,ampersand replacement=\&, column sep=0.5cm, row sep=0.5cm](m)
{
 \& $\cdots$ \& $H_{k+1}\left(\frac{N_I}{N}\right)$ \\
$H_k\left(\frac{N}{\partial_e N}\right)$ \& $H_k\left(\frac{N_I}{\partial_e N_I}\right)$ \& $H_k\left(\frac{N_I}{N}\right)$ \\
$H_{k-1}\left(\frac{N}{\partial_e N}\right)$ \& $\cdots$ \&  \\
};
\draw[->] (m-1-2) edge (m-1-3)
          (m-1-3) edge[curvedlink] (m-2-1)
          (m-2-1) edge (m-2-2) 
          (m-2-2) edge (m-2-3)
          (m-2-3) edge[curvedlink] (m-3-1)
          (m-3-1) edge (m-3-2);
\filldraw [black] (4.1075,-0.25) circle (0.01pt) node[anchor=west] {$d_{k-1,k}$};
\filldraw [black] (4.265,1.15) circle (0.01pt) node[anchor=west] {$d_{k,k+1}$};
\end{tikzpicture}
\end{center}
where the boundary maps $d_{k-1,k}$ constitute the entries of the matrix $\Delta_{i,i+1}=\delta_{i,i+1}$. Note this map is a flow defined boundary map, that is, it is not simply a topological map, but in its construction it takes account the dynamics of the neighborhoods of $M_{i+1},M_i$ and $M_I$, and how entagled those are. 

Since homology is homotopy invariant, and the Conley index is independent of the isolating neighborhoods, the quotient spaces within the homologies can be considered as their homotopy type, yielding
\begin{center}
\begin{tikzpicture}
\matrix[matrix of nodes,ampersand replacement=\&, column sep=0.5cm, row sep=0.5cm](m)
{
 \& $\cdots$ \& $H_{k+1}\left(h(M_{i})\right)$ \\
$H_k\left(h(M_{i-1})\right)$ \& $H_k\left(h(M_I)\right)$ \& $H_k\left(h(M_{i})\right)$ \\
$H_{k-1}\left(h(M_{i-1})\right)$ \& $\cdots$ \&  \\
};
\draw[->] (m-1-2) edge (m-1-3)
          (m-1-3) edge[curvedlink] (m-2-1)
          (m-2-1) edge (m-2-2) 
          (m-2-2) edge (m-2-3)
          (m-2-3) edge[curvedlink] (m-3-1)
          (m-3-1) edge (m-3-2);
\filldraw [black] (4.475,-0.3) circle (0.01pt) node[anchor=west] {$d_{k-1,k}$};
\filldraw [black] (4.665,0.9) circle (0.01pt) node[anchor=west] {$d_{k,k+1}$};
\end{tikzpicture}
\end{center}

Therefore, the boundary map is independent of the choice of an isolating neighborhood $N_I$ with $\partial_e N_I$, and its subisolating neighborhood $N$ with $\partial_e N=\partial_e N_I$.

In the case $M_{i+1},M_i$ consist of a single equilibria of Morse index $i+1,i$ respectively, then the nonzero elements in the above long exact sequence is as in Figure \ref{fig:dijconstr} since the Conley index of hyperbolic fixed points are spheres of appropriate dimension.
\begin{figure}[ht!]
\centering
\begin{center}
\begin{tikzpicture}
\matrix[matrix of nodes,ampersand replacement=\&, column sep=0.5cm, row sep=0.5cm](m)
{
$0$ \& $H_{i+1}\left(h(M_I)\right)$ \& $\mathbb{Z}$ \\
$\mathbb{Z}$ \& $H_{i}\left(h(M_I)\right)$ \& $0$ \&  \\
};
\draw[->] (m-1-1) edge (m-1-2)
          (m-1-2) edge (m-1-3)
          (m-1-3) edge[curvedlink] (m-2-1)
          (m-2-1) edge (m-2-2)
          (m-2-2) edge (m-2-3);

\draw[->] (m-1-1) edge (m-2-1)
          (m-1-2) edge (m-2-2)
          (m-1-3) edge (m-2-3);

\filldraw [black] (2.29,0.3) circle (0.01pt) node[anchor=west] {$d_{i,i+1}$};
\end{tikzpicture}
\end{center}
\caption{The construction of the entry $d_{i,i+1}$ within the connection matrix $\Delta$ for two adjacent equilibria of Morse index $i+1$ and $i$.\label{fig:dijconstr}}
\end{figure}

Now, we will use the diagram chasing from the zig-zag lemma \ref{thmzigzag} in Section \ref{sec:alg}, in order to obtain dynamical information from the topological construction of $d_{i,i+1}$. Indeed, $d_{i,i+1}$ maps a ball from the unstable manifold of the equilibrium in $M_{i+1}$ to the unstable manifold of the equilibrium in $M_{i}$ through a global passage within the intricacies of the invariant set $M_I$. 

Recall the proof of the zig-zag lemma. We consider a $[c]\in H_{i+1}(h(M_{i+1}))=\mathbb{Z}$. Therefore, there is a nontrivial linear combinations of $(i+1)$-cells $c\in C_{i+1}(h(M_{i+1}))$. Interpreting that $h(M_{i+1})$ is locally the unstable manifold of $M_{i+1}$, then 
\begin{equation}
    c\in W^u(M_{i+1}).    
\end{equation}

Note that $c\in [N]$ can not occur, since $c\in C_{i+1}$ is a cell of dimension $i+1$, and $[N]$ is a cell of dimension 0, i.e., a point in the quotient topology.

Due to the zig-zag lemma, there is a $[a]\in H_i(h(M_i))$ such that $a\in C_i(h(M_i))$ is a nontrivial linear combination of $i$-cells and $d_{i,i+1}([c])=[a]$. Similarly to the above argument, we obtain that 
\begin{equation}
    a\in W^u(M_i).
\end{equation}

Moreover, since $M_I$ is the maximal invariant set in $N_I$, then $T(t)a\in \partial_e N_I$ for some $t\in\mathbb{R}$. Therefore, the boundary map $d_{i,i+1}$ maps a representative of an $(i+1)$-cell in $W^u(M_{i+1})$ to an $i$-cell in $W^u(M_{i})$. We will track this using the diagram chasing in zig-zag, in order to understand the global passage for the unstable flow nearby $M_{i+1}$ to the unstable flow nearby $M_i$.

Following zig-zag, there is a $b\in C_{i+1}(h(M_I))$ such that $p_{i+1}(b)=c$. Note that either $b\in C_{i+1}(N_I\backslash M_I)$, $b\in C_{i+1}(M_I)$ or $b\in [\partial_e N_I]$. The third option can not occur, since it consists of a zero dimensional cell in the quotient topology.

Note that 
\begin{equation}
    b\in N_I\backslash N.    
\end{equation}

On the contrary, if $b\in N$, then $p_{i+1}(b)\in [N]$. But we know that $p_{i+1}(b)=c\not\in [N]$, since $c$ is a nontrivial cycle.

Note that the $(i+1)$-dimensional cell $b$ contains the cell $c$ of same dimension, since $p_{i+1}$ is surjective. Moreover, if part of the cell $b$ lies in $N$, then this is quotient out in the space $N_I\backslash N$. Therefore, $b$ can be seen as a ``neighborhood" of $c$. 

Again, following zig-zag, we obtain part of the global dynamical passage: 
\begin{equation}
    \partial^{M_I}_{i+1}(b)\in Ker(p_i)
    \subseteq C_i(h(M_I)).    
\end{equation}

Interpreting this equation: $\partial^{M_I}_{i+1}(b)$ is a $i$-dimensional cell that lies in $Ker(p_i)$, and hence $p_i(\partial^{M_I}_{i+1}(b))=[N]$. Therefore 
$\partial^{M_I}_{i+1}$ yields the global passage that maps $b\in N_I\backslash N$ into 
\begin{equation}
    \partial^{M_I}_{i+1}(b)\in N.
\end{equation}

The last diagram-chase from the zig-zag, tells us that there is an unique nontrivial $i$-dimensional cell $a\in C_i(M_i)$ such that
\begin{equation}\label{X}
    inc_i(a)=\partial^{M_I}_{i+1}(b).
\end{equation}

In other words, the cells $inc_i(a)$ and $\partial^{M_I}_{i+1}(b)$ are the same cell in $h(M_I)$, i.e., in $N_I\backslash \partial_e N_I$. Moreover, since $inc_i$ is injective, it just includes the $i$-dimensional cell $a\in N_I$ into the same cell $inc_i(a)\in N_I\backslash \partial_e N_I$. 

Lastly, recall that $a\in C_i(h(M_i))$. Similarly to the argument for $c$, this implies that $a\in W^u(M_i)$. Hence, since $\partial^{M_I}_{i+1}(b)$ lies in the same $i$-cell as $a$ due to $\eqref{X}$, they are both contained in the same $i$-dimensional ball $B^i$. Therefore, we can use the $\lambda$-lemma in order to track back such ball $B^i$ to a neighborhood of the hyperbolic fixed point $M_i$ through the semiflow $T(-t)B^i$ for some $t$. Thus, there will be a point in $T(-t)B^i$ that lies in the stable manifold $W^s(M_i)$ that can be traced back to the unstable manfold $W^u(M_{i+1})$ through the global passage.

\begin{figure}[ht]
\minipage{0.25\textwidth}\centering
\begin{tikzpicture}[scale=1]
    \draw [domain=0:6.28,variable=\t,smooth] plot ({sin(\t r)},{1.6*cos(\t r)}) node[anchor=south] {$N_I$}; 
    \draw [fill=lightgray,domain=1.45:4.85,variable=\t,smooth] plot ({sin(\t r)},{1.6*cos(\t r)}) ;   
    \draw [line width=0.6mm, domain=3.14:4.4,variable=\t,smooth] plot ({sin(\t r)},{1.6*cos(\t r)}) ;   
    \draw [line width=0.6mm,domain=1.88:3.14, variable=\t,smooth] plot ({sin(\t r)},{1.6*cos(\t r)}) node[anchor=north]{$\partial_e N_I=\partial_e N$};
    
    \filldraw [black] (0,-1) circle (1pt) node[anchor=west]{\tiny{$M_{i}$}};
    
    \filldraw [black] (0,1) circle (1pt) node[anchor=west] {\tiny{$M_{i+1}$}};    
    \draw [->,domain=-1:-0.5,variable=\t,smooth] plot ({0},{(-\t)}) node[right] {\tiny{$W^u$}};
      
    \filldraw [black] (0,0.75) circle (1pt) node[anchor=east] {\tiny{$c$}};    
    
    \draw (-1,0.2) -- (1,0.2);
    \filldraw [black] (-1.2,0.2) circle (0.01pt) node[anchor=north]{$N$};       
\end{tikzpicture}
\endminipage\hfill
\minipage{0.25\textwidth}\centering
\begin{tikzpicture}[scale=1]
    \draw [domain=0:6.28,variable=\t,smooth] plot ({sin(\t r)},{1.6*cos(\t r)}) node[anchor=south] {$N_I$}; 
    \draw [fill=lightgray,domain=1.45:4.85,variable=\t,smooth] plot ({sin(\t r)},{1.6*cos(\t r)}) ;   
    \draw [line width=0.6mm, domain=3.14:4.4,variable=\t,smooth] plot ({sin(\t r)},{1.6*cos(\t r)}) ;   
    \draw [line width=0.6mm,domain=1.88:3.14, variable=\t,smooth] plot ({sin(\t r)},{1.6*cos(\t r)}) node[anchor=north]{$\partial_e N_I=\partial_e N$};
    
    \filldraw [black] (0,-1) circle (1pt) node[anchor=west]{\tiny{$M_{i}$}};
    
    \filldraw [black] (0,1) circle (1pt) node[anchor=west] {\tiny{$M_{i+1}$}};    
    \draw [fill=lightgray,domain=0:6.28,variable=\t,smooth] plot ({0.1*sin(\t r)},{0.75+0.1*cos(\t r)}) node[anchor=west] {\tiny{$b$}};       
    \draw [->,domain=-1:-0.5,variable=\t,smooth] plot ({0},{(-\t)}) node[right] {\tiny{$W^u$}};
    \filldraw [black] (0,0.75) circle (1pt) node[anchor=east] {\tiny{$c$}};    
    
    \draw (-1,0.2) -- (1,0.2);
    \filldraw [black] (-1.2,0.2) circle (0.01pt) node[anchor=north]{$N$};           
\end{tikzpicture}
\endminipage\hfill
\minipage{0.25\textwidth}\centering
\begin{tikzpicture}[scale=1]
    \draw [domain=0:6.28,variable=\t,smooth] plot ({sin(\t r)},{1.6*cos(\t r)}) node[anchor=south] {$N_I$}; 
    \draw [fill=lightgray,domain=1.45:4.85,variable=\t,smooth] plot ({sin(\t r)},{1.6*cos(\t r)}) ;   
    \draw [line width=0.6mm, domain=3.14:4.4,variable=\t,smooth] plot ({sin(\t r)},{1.6*cos(\t r)}) ;   
    \draw [line width=0.6mm,domain=1.88:3.14, variable=\t,smooth] plot ({sin(\t r)},{1.6*cos(\t r)}) node[anchor=north]{$\partial_e N_I=\partial_e N$};
    
    \filldraw [black] (0,-1) circle (1pt) node[anchor=west]{\tiny{$M_{i}$}};
    
    \filldraw [black] (0,1) circle (1pt) node[anchor=west] {\tiny{$M_{i+1}$}};    
    \draw [fill=lightgray,domain=0:6.28,variable=\t,smooth] plot ({0.1*sin(\t r)},{0.75+0.1*cos(\t r)}) node[anchor=west] {\tiny{$b$}};   
    \draw [domain=-1:-0.5,variable=\t,smooth] plot ({0},{(-\t)}) node[right] {\tiny{$W^u$}};
    \draw [->,domain=-0.5:0.3,variable=\t,smooth] plot ({0},{(-\t)});
    
    \draw [domain=0:6.28,variable=\t,smooth] plot ({0.1*sin(\t r)},{-0.25+0.1*cos(\t r)}) node[anchor=west] {\tiny{$\partial_{i+1}^{M_I}(b)$}};   

    \draw (-1,0.2) -- (1,0.2);
    \filldraw [black] (-1.2,0.2) circle (0.01pt) node[anchor=north]{$N$};                
\end{tikzpicture}
\endminipage\hfill
\minipage{0.25\textwidth}\centering
\begin{tikzpicture}[scale=1]
    \draw [domain=0:6.28,variable=\t,smooth] plot ({sin(\t r)},{1.6*cos(\t r)}) node[anchor=south] {$N_I$}; 
    \draw [fill=lightgray,domain=1.45:4.85,variable=\t,smooth] plot ({sin(\t r)},{1.6*cos(\t r)}) ;   
    \draw [line width=0.6mm, domain=3.14:4.4,variable=\t,smooth] plot ({sin(\t r)},{1.6*cos(\t r)}) ;   
    \draw [line width=0.6mm,domain=1.88:3.14, variable=\t,smooth] plot ({sin(\t r)},{1.6*cos(\t r)}) node[anchor=north]{$\partial_e N_I=\partial_e N$};
    
    \filldraw [black] (0,-1) circle (1pt) node[anchor=west]{\tiny{$M_{i}$}};
    \draw [->,domain=1:1.5,variable=\t,smooth] plot ({0},{(-\t)});
    \filldraw [black] (0,-1.25) circle (1pt) node[anchor=east] {\tiny{$a$}};    
    \filldraw [black] (0,-1.25) circle (1pt) node[right] {\tiny{$W^u$}};
    
    \filldraw [black] (0,1) circle (1pt) node[anchor=west] {\tiny{$M_{i+1}$}};    
    \draw [fill=lightgray,domain=0:6.28,variable=\t,smooth] plot ({0.1*sin(\t r)},{0.75+0.1*cos(\t r)}) node[anchor=west] {\tiny{$b$}};   
    \draw [domain=-1:-0.5,variable=\t,smooth] plot ({0},{(-\t)}) node[right] {\tiny{$W^u$}};
    \draw [->,domain=-0.5:0.3,variable=\t,smooth] plot ({0},{(-\t)});
    
    \draw [domain=0:6.28,variable=\t,smooth] plot ({0.1*sin(\t r)},{-0.25+0.1*cos(\t r)}) node[anchor=west] {\tiny{$inc_i(a)$}};  
    
    \draw (-1,0.2) -- (1,0.2);
    \filldraw [black] (-1.2,0.2) circle (0.01pt) node[anchor=north]{$N$};           
\end{tikzpicture}
\endminipage
\captionof{figure}{Dynamical interpretation of algebraic proof of the zig-zag lemma.}
\end{figure}
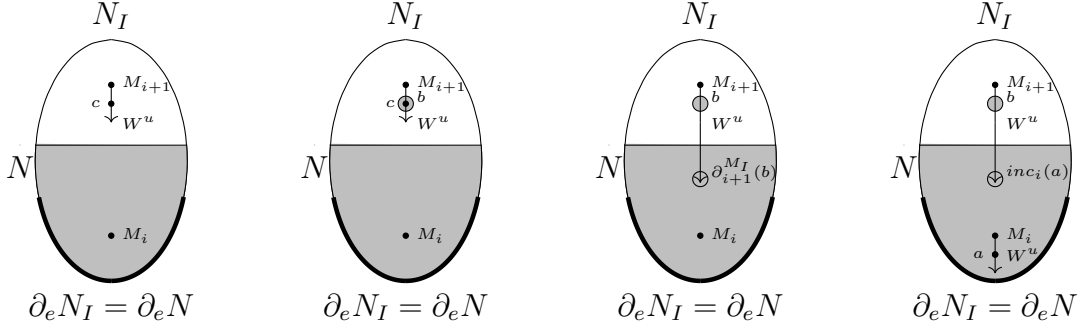

This implies that the zig-zag lemma yields a dynamical global passage through the abstract topological construction of the map $c\mapsto a$, that carries dynamical information by mapping a ball of dimension $i+1$ in $W^u(M_{i+1})$ into a ball of dimension $i$ in $W^u(M_{i})$. 

In particular, if the global pasnlsage map $\partial^{M_I}_{i+1}=0$, this implies that $\partial^{M_I}_{i+1}(b)\in [\partial_e N_I]$, and this means that the dynamics nearby $M_{i+1}$ that is described by the cell $c$ is carried over by the global passage until the exit set $\partial_e N_I$. And hence, this method is not able to track the point $\partial_e N_I(b)$ to be nearby $M_{i+1}$. Hence, the connection matrix is inconclusive in order to detect heteroclinics - even though there might be some.

On the other hand, if $\partial_e N_I\neq 0$, there are points in the unstable manifold $W^u(M_{i+1})$ lying in the cell $c$ such that are mapped through the global passage to the cell $a$ that contains points of $W^u(M_i)$. Then, using the $\lambda$-lemma argument above, we can find a point in the stable manifold $W^s(M_i)$ that can be tracked back to $W^u(M_{i+1})$, yielding a heteroclinic connection.

After this dynamical interpretation of algebraic facts, we finish this section with two remaks.

First, if there is no connection between $M_{i+1}$ and $M_i$, that is $M_I=\{M_{i+1},M_i\}$, then $h(M_I)=h(M_{i+1})\vee h(M_{i})$, and hence $H\left(h(M_I)\right) = H\left(h(M_{i+1})\right)\oplus H\left(h(M_{i+1})\right)$. Algebraically, this means that the sequence splits. And hence, the middle vertical map has to be zero, since it comes from the chain complex level maps. This yields that $d_{i,i+1}=0$. The counter-positive of the statement in this paragraph, is a proof of Theorem \ref{thmdetection} for this particular case.

Second, if we can compute $h(M_I),h(M_{i+1})$ and $h(M_{i})$, then one can show the existence of a connecting orbit between $M_{i+1}$ and $M_{i}$ within $M_I$. In fact, suppose towards a contradiction that there is no heteroclinic connection between $M_{i+1}$ and $M_i$, that is, $M_I=M_{i+1} \cup M_{i}$. Then, $h(M_I)=h(M_{i+1})\vee h(M_{i})$. Therefore, if we know  $h(M_I)\neq h(M_{i+1})\vee h(M_{i})$, then we obtain a contradiction and there is a connecting heteroclinic.

\section{Applications}\label{sec:app}

\subsection{Delay equation - periodic orbit connections}

In this section, we address the connection problem for certain differential equations with delayed. We follow Mischaikov \cite{Mischaikow87} in order to compute the connection matrix $\Delta_{i,j}$ using the knowledge of the Conley index of Morse sets and the properties $1$ to $4$ from Theorem \ref{ThmDelta} above. Hence, we conclude the existence of connecting orbits in $C_{i,j}$.

Consider the delay differential equation
\begin{equation}
    \dot{u}=f(u(t),u(t-1))
\end{equation}
with initial condition $u_0\in X:=C^0([-1,0])$.

It is known that solutions $u(t)$ generate a semiflow in $X$. Moreover, under appropriate growth conditions on $f$, there exists a global attractor $\mathcal{A}\subseteq X$, which is the maximal compact invariant set. See Hale and Lunel \cite{HaleLunel93}.

It was shown by Mallet-Paret \cite{MalletParet88} that there is a Morse decompotision of the global attractor, due to nodal properties. Instead of prescribing $f$, we prescribe the invariant set for the dynamics. For our example, a nonlinearity $f$ can be chosen, for instance, as one that has gone under consecutive Hopf bifurcations. An example of such bifurcations is the Wright's equation, when $f(u(t),u(t-1)):=-\alpha u(t-1)[1+u(t)]$ which has a supercritical Hopf bifurcatoin at $\alpha_n:=\pi/2+2\pi n$ for $n\in\mathbb{N}_0$. See 
and 

Suppose that the global attractor $\mathcal{A}$ has a Morse decomposition with Morse sets given by $M_\iota:=\{u_\iota(t) | u_\iota(t) = u_\iota(t+p_\iota)\}$ for $\iota=0,1$, and $M_2:=\{u\equiv 0\}$ such that the index for each Morse set is related to its unstable dimension: there are 4 positive eigenvalues for the hyperbolic fixed point $u\equiv 0$; and the linearization of the Poincaré map $P_\iota$ for the periodic orbits $u_\iota(t)$ has exactly $2\iota$ positive eigenvalues. Therefore, as we computed in the previous Section, the Conley index of the hyperbolic fixed point and the periodic orbits are \eqref{HypHomindex} and \eqref{PTHomindex}, respectrively:
\begin{equation}
    H_k(h(M_\iota))=
    \begin{cases}
        \mathbb{Z}_2  & \text{ if } k=2\iota,2\iota+1 \\
        0 & \text{ if } k\neq 2\iota,2\iota+1
    \end{cases}    
    \text{ , }
    H_k(h(M_2))=
    \begin{cases}
        \mathbb{Z}_2  & \text{ if } k=4 \\
        0 & \text{ if } k\neq 4.
    \end{cases}    
\end{equation}

Therefore, the connection matrix $\Delta: \oplus_{i=0}^2 H(h(M_i))\to \oplus_{i=0}^2 H(h(M_i))$. Due to the property 1 in Theorem \ref{ThmDelta}, the lower triangle is zero. Hence, this map can be represented as the following matrix form
\begin{equation}
  \Delta=
  \begin{blockarray}{*{3}{c} l}
    \begin{block}{*{3}{>{$\footnotesize}c<{$}} l}
      $H(h(M_0))$ & $H(h(M_1))$ & $H(h(M_2))$ \\
    \end{block}
    \begin{block}{[*{3}{c}]>{$\footnotesize}l<{$}}
      0 & \Delta_{01} & \Delta_{02} \bigstrut[t] & $H(h(M_0))$ \\
      0 & 0 & \Delta_{12} & $H(h(M_1))$ \\
      0 & 0 & 0 & $H(h(M_2))$. \\
    \end{block}
  \end{blockarray}
\end{equation}

Note that since each Morse set has only one connected component, the entries $\Delta_{ij}$ of the matrix $\Delta$ consist of a single entry, denoted by $\delta_{ij}$, which has entries $d_{mn}$ as below
\begin{equation}\label{DeltaMatrixDELAY}
  \Delta=
  \begin{blockarray}{*{5}{c} l}
    \begin{block}{*{5}{>{$\footnotesize}c<{$}} l}
      $H_0(h(M_0))$ & $H_1(h(M_0))$ & $H_2(h(M_1))$ & $H_3(h(M_{1}))$ & $H_4(h(M_2))$  \\
    \end{block}
    \begin{block}{[*{5}{c}]>{$\footnotesize}l<{$}}
      0 & 0 & d_{02} & d_{03} & d_{04} \bigstrut[t] & $H_0(h(M_0))$ \\
      0 & 0 & d_{12} & d_{13} & d_{14} & $H_1(h(M_0))$ \\
      0 & 0 & 0 & 0 & d_{24} & $H_2(h(M_1))$ \\
      0 & 0 & 0 & 0 & d_{34} & $H_3(h(M_1))$ \\
      0 & 0 & 0 & 0 & 0 & $H_4(h(M_2))$ \\
    \end{block}
  \end{blockarray}
\end{equation}

Property 3 in Theorem \ref{ThmDelta} implies that $d_{04}=d_{14}=d_{24}=d_{03}=d_{13}=d_{02}=0$, since $d_{ij}$ is a degree $-1$ map and maps $k^{th}$ to $(k-1)^{th}$ level homology.

The remaining possible nonzero entries are $d_{12}$ and $d_{34}$. Note that property 2 does not tell us anything new, since $\Delta^2=0$ for any value of $d_{12}$ and $d_{34}$.

We now examine property 4 for $I=\{0,1,2\}$, that is, $M_I=\mathcal{A}$ and $\Delta_I=\Delta$. Since the Conley index of the global attractor is known, $\dim H(h(\Sigma))=1$ using \eqref{ConleyAtt}, then property 4 implies that $\dim (\text{Ker}\Delta / \text{Im} \Delta)=1$. In other words, 
\begin{equation}\label{randomshit}
    \dim(\text{Ker}\Delta) - \dim(\text{Im} \Delta)=1.
\end{equation}

Since most $d_{ij}$ are zero within the matrix in \eqref{DeltaMatrixDELAY}, except possibly $d_{12}$ and $d_{34}$, we obtain that there are at least three null vertical columns. Therefore, $\dim(\text{Ker}\Delta)\geq 3$. Together with \eqref{randomshit}, it implies that $\dim(\text{Im}\Delta)\geq 2$, and consequently $\Delta$ must have at least rank 2. Finally, $d_{12},d_{34}\neq 0$ and Theorem \ref{thmdetection} implies that there are heteroclinic connections from $M_2$ to $M_1$, and from $M_1$ to $M_0$.

\begin{figure}[ht]\centering
    \begin{tikzpicture}[scale=1]

    \filldraw [black] (0,4) circle (2pt) node[anchor=south] {$u\equiv 0$}; 
    
    \draw[thick,->] (0,4) -- (0,3.1);

    \draw (-1,2.5) arc (180:360:1cm and 0.5cm);
    \draw (-1,2.5) arc (180:0:1cm and 0.5cm);
    \draw[->] (-1,2.5) -- (-1,2.49) node[anchor=east]{$M_1$};
    
    \draw[thick,->] (0,2) -- (0,1.1);

    \draw (-1,0.5) arc (180:360:1cm and 0.5cm);
    \draw (-1,0.5) arc (180:0:1cm and 0.5cm);
    \draw[->] (-1,0.5) -- (-1,0.49) node[anchor=east]{$M_0$};
    
    \end{tikzpicture}
\caption{The global attractor $\mathcal{A}$ for a delay equation with heteroclinics between the equilibrium $u\equiv 0$ and the periodic orbit $M_1$, and also between the periodic orbits $M_1$ and $M_0$.}
\label{attDELAY}
\end{figure}
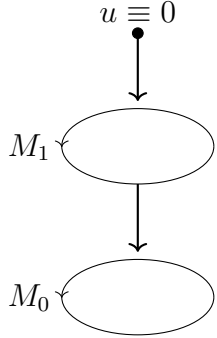

\subsection{Symmetric bistable gradient-like systems - Chafee-Infante type connections}

We follow some ideas within Mischaikow \cite{Mischaikow95}, which proves more than is presented here: we show the existence of heteroclinics between equilibria of bistable gradient-like systems (yet to be defined), whereas Mischaikow shows those systems are topologically semiconjugate to the flow of an ODE in the $n$-dimensional ball, in particular, consisting of the heteroclinics we compute here.

We consider the semiflow $T(t)$ as in \eqref{semiflow}, which depends also on a parameter $\lambda\in \mathbb{R}_+$ continuously. From now on, we denote it by $T_\lambda (t)$. We suppose that $T_\lambda (t)$ has a compact global attractor $\mathcal{A}_\lambda$, which is also parameter dependent.

The semiflow $T_\lambda(t)$ is \emph{symmetric} if for any $u_0\in X$, then $T(t)(-u_0)=-T(t)u_0$. This assumption is not made by Mischaikow, nevertheless this symmetry is inherited by the Chafee-Infante problem and will be used in the proof given here.

The semiflow $T_\lambda(t)$ is \emph{bistable} if there are two stable hyperbolic equilibria $e^-_0,e^+_0$ for all $\lambda$, and a third equilibrium $e_N$ that is hyperbolic for all $\lambda$, except for a sequence $\lambda_n$ with $n\in \mathbb{N}_0$ such that $\lambda_n<\lambda_{n+1}$ and $\lambda_n\to\infty$. Moreover, at the parameter value $\lambda_n$, the equilibrium $e_N$ undergoes a supercritical pitchfork bifurcation, yielding two new hyperbolic equilibria $e^-_n,e^+_n$. In total, the set of equilibria consists of $2n+1$ hyperbolic elements given by $\mathcal{E}_\lambda=\{e_0^\pm,e_1^\pm,...,e^\pm_{n},e_N\}$ for $\lambda\in (\lambda_n,\lambda_{n+1})$. Note that hyperbolicity means that no other bifurcations can occur. 

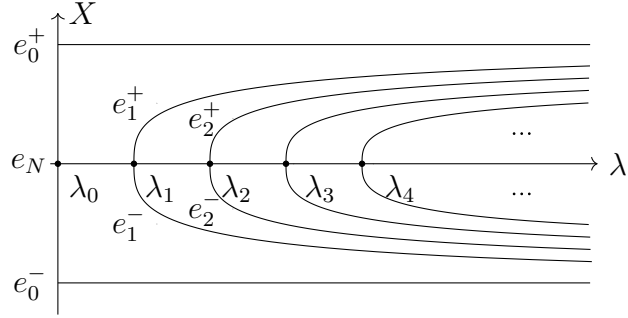
\begin{figure}[ht!]
\centering
\begin{tikzpicture}[scale=1]
        \draw[->] (-0.1,0) -- (7.1,0) node[right] {$\lambda$};
        \draw[->] (0,-2) -- (0,2) node[right] {$X$};
        
        \draw[-] (7,1.58) -- (0,1.58) node[left] {$e_0^+$};        
        \draw[-] (7,-1.58) -- (0,-1.58) node[left] {$e_0^-$};
        \filldraw [black] (0,0) circle (1pt) node[anchor=north west]{$\lambda_0$};
        \filldraw [black] (0,0) circle (0.01pt) node[anchor=east]{$e_N$};        
        \draw[domain=-1.298:1.298,samples=100] plot ({1+\x^2*tan((\x)r)},\x);
        \filldraw [black] (1.3,0.8) circle (0.01pt) node[anchor=east]{$e^+_1$};
        \filldraw [black] (1.3,-0.8) circle (0.01pt) node[anchor=east]{$e^-_1$};
        \filldraw [black] (1,0) circle (1pt) node[anchor=north west]{$\lambda_1$};

        \draw[domain=-1.2625:1.2625,samples=100] plot ({2+\x^2*tan((\x)r)}, 0.9*\x);
        \filldraw [black] (2.3,0.6) circle (0.01pt) node[anchor=east]{$e^+_2$};
        \filldraw [black] (2.3,-0.6) circle (0.01pt) node[anchor=east]{$e^-_2$};
        \filldraw [black] (2,0) circle (1pt) node[anchor=north west]{$\lambda_2$};

        \draw[domain=-1.2165:1.2165,samples=100] plot ({3+\x^2*tan((\x)r)}, 0.8*\x);
        \filldraw [black] (3,0) circle (1pt); 
        \filldraw [black] (3.1,0) circle (0.01pt) node[anchor=north west]{$\lambda_3$};

        \draw[domain=-1.152:1.152,samples=100] plot ({4+\x^2*tan((\x)r)}, 0.7*\x);
        \filldraw [black] (4,0) circle (1pt);
        \filldraw [black] (4.15,0) circle (0.01pt) node[anchor=north west]{$\lambda_4$};

        \filldraw [black] (6,0.4) circle (0.3pt);
        \filldraw [black] (6.1,0.4) circle (0.3pt);
        \filldraw [black] (6.2,0.4) circle (0.3pt);      

        \filldraw [black] (6,-0.4) circle (0.3pt);
        \filldraw [black] (6.1,-0.4) circle (0.3pt);
        \filldraw [black] (6.2,-0.4) circle (0.3pt);     

\end{tikzpicture}
\caption{This picture depicts the two stable hyperbolic equilibria $e_0^\pm$, and the supercritical pitchfork bifurcation of the equilibrium $e_N$ at each parameter value $\lambda_n$, yielding the set of all equilbiria $\mathcal{E}_\lambda=\{e_0^\pm,e_1^\pm,...,e^\pm_{n},e_N\}$ for $\lambda\in (\lambda_n,\lambda_{n+1})$.}
\end{figure}

The semiflow $T_\lambda(t)$ is \emph{bistable gradient-like} if semiflow $T_\lambda(t)$ is bistable, and there is a Lyapunov function $L:\mathcal{A}\to\mathbb{R}$ that is decreasing along trajectories, except at equilibria $\mathcal{E}_\lambda$ such that $L(e^-_i)=L(e^+_i)$ for all $i=0,...,n$.

Several examples of bistable systems are listed in \cite{Mischaikow95}: the Chafee-Infante heat equation, see \cite{ChafeeInfante74}; a damped wave equation \cite{Hale88}; a one-dimensional beam with soft loading, see \cite{HattoriMischaikow91}; a FitzHugh-Nagumo equation \cite{ConleySmoller86}; a Cahn-Hilliard putting together \cite{Temam88} and \cite{BatesFife90}; a phase-field equations \cite{BatesZheng92}.

The prototype we will focus is the Chafee-Infante heat equation
\begin{equation}\label{CI}
    u_t=u_{xx}+\lambda u(1-u^2)
\end{equation}
for $x\in (0,1)$ with Neumann boundary conditions. 

The results in this section are valid if we replace $u(1-u^2)$ with more general nonlinearity $f(u)\in C^2$ that have a similar shape: $f(0)=0$, $f'(0)=1$, $\limsup_{|u|\to\infty} f(u)/u \leq 0$ and $uf''(u)<0$ for all $x\in \mathbb{R}^2\backslash \{0\}$. Those yield a sequence of pitchfork bifurcations and have a Lyapunov function.

The Chafee-Infante problem has two stable hyperbolic equilibria $e^\pm_0=\pm 1$, and an equilibrium $e_N\equiv 0$ that undergoes a supercritical pitchfork bifurcation at $\lambda_n=n^2$. See \cite{ChafeeInfante74}. Moreover, it has a Lyapunov function given by 
\begin{equation}
    L_\lambda(u,u_x)=\int_0^1 \frac{|u_x|^2}{2}-\lambda \left( \frac{u^2}{2}-\frac{u^4}{4} \right)dx \text{\hspace{0.5cm} such that \hspace{0.5cm}} \frac{dL_\lambda}{dt}=-\int_0^1 |u_t|^2dx.
\end{equation}

Note that the symmetry assumption of the semiflow is manifested in equation \eqref{CI} as the symmetry $u\to -u$, namely $u$ is a solution if, and only if $-u$ is as well. Therefore, the equilibria that arise from the pitchfork bifurcation are also symmetric: $e^+_i=-e^-_i$, and consequently $L(e^+_i)=L(e^-_i)$ for all $i=0,...,n$.

The Lyapunov function implies the LaSalle invariance principle holds, see \cite{Henry81}, and hence the $\alpha,\omega$-limits consist of a single equilibrium. Following the notation in Section \ref{sec:matrix}, there is a Morse decomposition for $\lambda\in (\lambda_n,\lambda_{n+1})$ given by $\mathcal{M}:=\{M_0,M_1,...,M_n,M_N\}$ where $M_i:=\{e_i^-\}\cup\{e_i^+\}$ for $i=0,...,n$, and $M_N:=\{e_N\}$.

\begin{prop}
If $T_\lambda(t)$ is a symmetric bistable gradient-like system for $\lambda\in (\lambda_n,\lambda_{n+1})$, then connection matrix of the form \eqref{DeltaMatrix} is given by 
\begin{equation}\label{DeltaMatrixCI}
  \Delta=
  \begin{blockarray}{*{5}{c} l}
    \begin{block}{*{5}{>{$\footnotesize}c<{$}} l}
      $H(h(M_0))$ & $H(h(M_1))$ & $H(h(M_2))$ & $\dots$ & $H(h(M_N))$  \\
    \end{block}
    \begin{block}{[*{5}{c}]>{$\footnotesize}l<{$}}
      0 & \Delta_{01} & 0 & \dots & 0 \bigstrut[t] & $H(h(M_0))$ \\
      \vdots & \ddots & \ddots & \ddots & \vdots & $\vdots$ \\
      0 & 0 & 0 & \ddots & 0 & $H(h(M_{n-1}))$ \\
      0 & 0 & 0 & \ddots & \Delta_{nN} & $H(h(M_n))$\\
      0 & 0 & 0 & \dots & 0 & $H(h(M_N))$ \\
    \end{block}
  \end{blockarray}
\end{equation}
where the upper diagonal submatrices $\Delta_{i,i+1}:H_{i+1}(h(M_{i+1}))\to H_{i}(h(M_{i}))$
for $i=0,...,n-1$ are only nonzero entries,
\begin{equation}\label{DeltaijCI2}
  \Delta_{i,i+1}=
  \begin{blockarray}{*{2}{c} l}
    \begin{block}{*{2}{>{$\footnotesize}c<{$}} l}
      $H_{i+1}(h(e^+_{i+1}))$ & $H_{i+1}(h(e^-_{i+1}))$\\
    \end{block}
    \begin{block}{[*{2}{c}]>{$\footnotesize}l<{$}}
      d_{++} & d_{+-} \bigstrut[t] & $H_{i}(h(e^+_{i}))$ \\
      d_{-+}  & d_{--}  & $H_{i}(h(e^-_{i}))$ \\
    \end{block}
  \end{blockarray}
\end{equation}
with all $d_{\pm\pm}\neq 0$, and $\Delta_{n,N}:H_{N}(h(M_{N}))\to H_{n}(h(M_{n}))$ is
\begin{equation}\label{DeltaijCI3}
  \Delta_{n,N}=
  \begin{blockarray}{*{1}{c} l}
    \begin{block}{*{1}{>{$\footnotesize}c<{$}} l}
      $H_{N}(h(e_{N}))$\\
    \end{block}
    \begin{block}{[*{1}{c}]>{$\footnotesize}l<{$}}
      d_{+N} \bigstrut[t] & $H_{N-1}(h(e^+_{n}))$ \\
      d_{-N} \bigstrut[t] & $H_{n}(h(e^-_{n}))$ \\      
    \end{block}
  \end{blockarray}
\end{equation}
for $d_{\pm N}\neq 0$.
\end{prop} 

From the above proposition, we can infer the following.

\begin{cor}
Consider a symmetric bistable gradient-like system $T_\lambda (t)$ for $\lambda\in (\lambda_n,\lambda_{n+1})$ with equilibria $\mathcal{E}_\lambda=\{e_0^\pm,e_1^\pm,...,e^\pm_{n},e_N\}$. Then, there are heteroclinic connections from the equilibria $e^\pm_{i+1}$ to $e^\pm_{i}$ for all $i=0,...,n$, and from $e_N$ to $e^\pm_{n}$.
\end{cor}
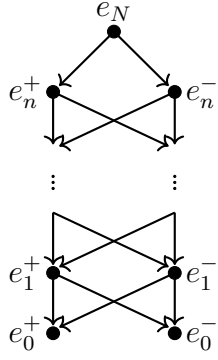
\begin{figure}[ht]\centering
\begin{tikzpicture}[scale=0.8]
\filldraw [black] (0,0) circle (3pt) node[anchor=south]{$e_{N}$};
\filldraw [black] (-1,-1) circle (3pt) node[anchor=east]{$e^+_{n}$};
\filldraw [black] (1,-1) circle (3pt) node[anchor=west]{$e^-_{n}$};

\filldraw [black] (-1,-2.4) circle (0.5pt);
\filldraw [black] (-1,-2.5) circle (0.5pt);
\filldraw [black] (-1,-2.6) circle (0.5pt);
\filldraw [black] (1,-2.4) circle (0.5pt);
\filldraw [black] (1,-2.5) circle (0.5pt);
\filldraw [black] (1,-2.6) circle (0.5pt);

\filldraw [black] (-1,-4) circle (3pt) node[anchor=east]{$e^+_{1}$};
\filldraw [black] (1,-4) circle (3pt) node[anchor=west]{$e^-_{1}$};
\filldraw [black] (-1,-5) circle (3pt) node[anchor=east]{$e^+_{0}$};
\filldraw [black] (1,-5) circle (3pt) node[anchor=west]{$e^-_{0}$};

\draw[thick,->] (0,0) -- (-0.9,-0.9);
\draw[thick,->] (0,0) -- (0.9,-0.9);

\draw[thick,->] (-1,-1) -- (0.9,-1.9);
\draw[thick,->] (1,-1) -- (-0.9,-1.9);
\draw[thick,->] (-1,-1) -- (-1,-1.87);
\draw[thick,->] (1,-1) -- (1,-1.87);

\draw[thick,->] (-1,-3) -- (0.9,-3.9);
\draw[thick,->] (1,-3) -- (-0.9,-3.9);
\draw[thick,->] (-1,-3) -- (-1,-3.87);
\draw[thick,->] (1,-3) -- (1,-3.87);

\draw[thick,->] (-1,-4) -- (0.9,-4.9);
\draw[thick,->] (1,-4) -- (-0.9,-4.9);
\draw[thick,->] (-1,-4) -- (-1,-4.87);
\draw[thick,->] (1,-4) -- (1,-4.87);

\end{tikzpicture}
\caption{Depiction of bistable gradient-like global attractor $\mathcal{A}_\lambda$ for $\lambda\in (\lambda_n,\lambda_{n+1})$} \label{FIGCOR}
\end{figure}

The proof is divided in two parts: we calculate the Conley index of each Morse set, and then we compute the connection matrix using the properties from Theorem \ref{thmdetection}.

Firstly, note that each hyperbolic equilibria $e_i^\pm$ with $i=0,...,n,N$ has Morse index $i$ for $\lambda\in (\lambda_n,\lambda_{n+1})$, and hence we can compute its homological Conley index, as in \eqref{HypHomindex}, yielding
\begin{equation}\label{CIindex1}
    H_k(h(M_i))=H_k(\mathbb{S}^i \vee \mathbb{S}^i)=
    \begin{cases}
        \mathbb{Z}_2\oplus \mathbb{Z}_2  & \text{ if } k=i\\
        0 & \text{ if } k\neq i
    \end{cases}
\end{equation} 
for $i=0,...,n$, and 
\begin{equation}\label{CIindex2}
    H_k(h(M_N))=H_k(\mathbb{S}^{n+1})=
    \begin{cases}
        \mathbb{Z}_2  & \text{ if } k=n+1\\
        0 & \text{ if } k\neq n+1.
    \end{cases}
\end{equation} 

Indeed, we show this by induction on $n$. For $n=0$, no bifurcation occurred, and there are two stable equilibria $e^\pm_0$, and one unstable $e_N$ of Morse index 1. The induction step follows by a close analysis of the pitchfork bifurcation. Suppose \eqref{CIindex1} and \eqref{CIindex2} are valid for $\lambda\in (\lambda_n,\lambda_{n+1})$ and we will show it  remains valid for $\lambda\in (\lambda_{n+1},\lambda_{n+2})$. Due to hyperbolicity of the equilibria $e^\pm_i$ for $i=0,...,n$, the induction hypothesis implies \eqref{CIindex1} remains valid for those $i$. 
 As $\lambda$ crosses $\lambda_{n+1}$, two more new equilibria $e^\pm_{n+1}$ with Morse index $n+1$ appear due to the pitchfork, whereas $e_N$ increases its Morse index by one. Therefore \eqref{CIindex1} is also valid for $i=n+1$, and \eqref{CIindex2} is valid by switching $n+1$ to $n+2$. This proves the desired claim.
 
Secondly, we compute the connection matrix using the properties from Theorem \ref{thmdetection}.

Since $\Delta$ is strict upper triangular, then $\Delta_{ij}=0$ for $j\leq i$. Moreover, since it is a degree $-1$ map, and each Morse set $M_i$ only has nontrivial homology in the $i^{th}$ level, as in \eqref{CIindex1}, then $\Delta_{ij}=0$ for all $j>i+1$. Therefore, we obtain \eqref{DeltaMatrixCI}, since the only remaining possibly nonzero entries are the upper-diagonal ones given by $\Delta_{i,i+1}$ as 
\begin{equation}\label{DeltaijCIproof}
  \Delta_{i,i+1}=
  \begin{blockarray}{*{2}{c} l}
    \begin{block}{*{2}{>{$\footnotesize}c<{$}} l}
      $H_{i+1}(h(e^+_{i+1}))$ & $H_{i+1}(h(e^-_{i+1}))$\\
    \end{block}
    \begin{block}{[*{2}{c}]>{$\footnotesize}l<{$}}
      d_{++} & d_{+-} \bigstrut[t] & $H_{i}(h(e^+_{i}))$ \\
      d_{-+} & d_{--} & $H_{i}(h(e^-_{i}))$ \\
    \end{block}
  \end{blockarray}
\end{equation}
and 
\begin{equation}\label{DeltaijCI4}
  \Delta_{n,N}=
  \begin{blockarray}{*{1}{c} l}
    \begin{block}{*{1}{>{$\footnotesize}c<{$}} l}
      $H_{N}(h(e_{N}))$\\
    \end{block}
    \begin{block}{[*{1}{c}]>{$\footnotesize}l<{$}}
      d_{+N} \bigstrut[t] & $H_{n}(h(e^+_{n}))$ \\
      d_{-N} \bigstrut[t] & $H_{n}(h(e^-_{n}))$. \\      
    \end{block}
  \end{blockarray}
\end{equation}

We proceed to show that the entries of the above matrices \eqref{DeltaijCIproof} and \eqref{DeltaijCI4} are nonzero by induction on $n$.  

For the induction basis, $n=0$, and hence $\lambda\in (\lambda_0,\lambda_1)$. In this case, the global attractor $\mathcal{A}_\lambda$ has Morse decomposition $\mathcal{M}:=\{M_0,M_N\}$ where $M_0:=\{e_0^-\}\cup\{e_0^+\}$ and $M_N:=\{e_N\}$. Hence the connection matrix $\Delta$ is given by
\begin{equation}\label{matrix1}
  \Delta=
  \begin{blockarray}{*{3}{c} l}
    \begin{block}{*{3}{>{$\footnotesize}c<{$}} l}
      $H_{0}(h(e^+_{0}))$ & $H_{0}(h(e^-_{0}))$ & $H_1(h(e_N))$ \\
    \end{block}
    \begin{block}{[*{3}{c}]>{$\footnotesize}l<{$}}
      0 & 0 & d_{+N} \bigstrut[t] & $H_{0}(h(e^+_{0}))$ \\
      0 & 0 & d_{-N} & $H_{0}(h(e^-_{0}))$ \\
      0 & 0 & 0 & $H_1(h(e_N))$. \\
    \end{block}
  \end{blockarray}
\end{equation}

We now prove that $d_{+N},d_{-N}\neq 0$. Choosing the interval $I:=\{0,N\}$, we have that $M_I=\mathcal{A}$ and $\Delta_I=\Delta$. Then the knowledge of the homology of the attractor \eqref{ConleyAtt}, and property 4 in Theorem \ref{thmdetection} implies
\begin{equation}\label{eq1}
    1=\dim(H(h(\mathcal{A})))=\dim(Ker \Delta)-\dim (Im(\Delta)).
\end{equation}

Since there are two zero columns in \eqref{matrix1}, we know $\dim(Ker \Delta)\geq 2$. Note it cannot be 3: in this case, $d_{+N},d_{-N}=0$ and hence \eqref{eq1} implies $1=3-0$, a contradiction. Therefore, $\dim(Ker \Delta)=2$. Also, $\dim (Im(\Delta))=1$. This implies that $d_{+N}$ and $d_{-N}$ cannot be both zero. 

Now we show that $d_{+N}$ is a multiple of $d_{-N}$. Recall how the entry $d_{+N}$ is constructed, finishing in Figure \eqref{fig:dijconstr}. For such construction, it is enough if we choose $M_{i+1}=\{e_N\}$ and $M_{i}=\{e^+_0\}$ with isolated invariant set $\Sigma:=\{e_N,e^+_0\}\cup \mathcal{C}^+_{N,0}$. Due to the symmetry of the semiflow, we can consider the same construction with $-M_{i+1}=M_{i+1}$, $-M^+_{i}=\{e^-_0\}$ and $-\Sigma=\{e_N,e^-_0\}\cup \mathcal{C}^-_{N,0}$, which will construct $d_{-N}$ in the exactly same way. Therefore, $d_{+N}$ is a multiple of $d_{-N}$.

The last two paragraphs conclude that both $d_{+N},d_{-N}\neq 0$.






Before we proceed with the induction step, we understand what occurs at the pitchfork bifurcation as the parameter $\lambda\in (\lambda_0,\lambda_1)$ crosses the value $\lambda_1$, and what the connection matrix becomes when  $\lambda\in (\lambda_1,\lambda_2)$. Then, the induction step will mimic this proof.

For $\lambda\in(\lambda_1,\lambda_2)$, the global attractor $\mathcal{A}_\lambda$ has Morse decomposition  $\mathcal{M}:=\{M_0,M_1,M_N\}$ where $M_0:=\{e_0^-\}\cup\{e_0^+\}$,$M_1:=\{e_1^-\}\cup\{e_1^+\}$ and $M_N:=\{e_N\}$. Therefore, the connection matrix is 
\begin{equation}\label{matrix2}
  \Delta=
  \begin{blockarray}{*{5}{c} l}
    \begin{block}{*{5}{>{$\footnotesize}c<{$}} l}
      $H_{0}(h(e^+_{0}))$ & $H_{0}(h(e^-_{0}))$ & $H_{1}(h(e^+_{1}))$ & $H_{1}(h(e^-_{1}))$ & $H_2(h(e_N))$  \\
    \end{block}
    \begin{block}{[*{5}{c}]>{$\footnotesize}l<{$}}
      0 & 0 & d_{++} & d_{+-} & 0 \bigstrut[t] & $H_{0}(h(e^+_{0}))$ \\
      0 & 0 & d_{-+} & d_{--} & 0 & $H_{0}(h(e^-_{0}))$ \\
      0 & 0 & 0 & 0 & d_{+N} & $H_{1}(h(e^+_{1}))$ \\
      0 & 0 & 0 & 0 & d_{-N} & $H_{1}(h(e^-_{1}))$ \\
      0 & 0 & 0 & 0 & 0 & $H_2(h(e_N))$ \\
    \end{block}
  \end{blockarray}
\end{equation}
where the entries above are zero because the connection matrix is strictly upper triangular, and a degree $-1$ map. Following \eqref{DeltaijCIproof} and \eqref{DeltaijCI4}, we denote the possibly nonzero square block by $\Delta_{0,1}:H_{1}(h(e^+_{1}))\oplus H_{1}(h(e^-_{1}))\to H_{0}(h(e^+_{0}))\oplus H_{0}(h(e^-_{0}))$, and the possibly nonzero two rowed last column $\Delta_{1,2}:H_{2}(h(e_{N}))\to H_{1}(h(e^+_{1}))\oplus H_{1}(h(e^-_{1}))$.  

Using that $\Delta$ is a boundary map $\Delta^2=0$, as in property 2 of Theorem \ref{thmdetection}, we obtain that $d_{++}=d_{+N}d_{+-}$ and $d_{-+}=d_{-N}d_{--}$. This means that the third and fourth column of \eqref{matrix2} are linearly dependent.

Since first two columns in \eqref{matrix2} are zero, and the third and fourth columns are linearly dependent, the dimension of the image (which is the span of the column vectors) of $\Delta$ is at most 2, $\dim Im (\Delta)\leq 2$, whereas the dimension of the kernel (the complement of the image, due to the rank-nullity theorem) is at least 3, $\dim Ker (\Delta)\geq 3$.

On the other hand, choose the interval $I:=\{0,1,N\}$, so that $M_I=\mathcal{A}$ and $\Delta_I=\Delta$. Then the knowledge of the homology of the attractor \eqref{ConleyAtt}, and property 4 in Theorem \ref{thmdetection} implies the same relation \eqref{eq1}. Using that the dimension of the kernel is at least 3 from the last paragraph, we obtain that the dimension of the image is at least 2, $\dim Im (\Delta)\geq 2$. Therefore, we obtain $\dim Im (\Delta)=2$. Similarly, $\dim Ker (\Delta)=3$.

Note that $d_{+N},d_{-N}$ cannot be both zero: otherwise, the dimension of the kernel will be at least 4, a contradiction. Moreover, using the symmetry of the semiflow as above, $d_{+N}$ and $d_{-N}$ are constructed in the same way, and hence $d_{+N}=0$ if, and only if $d_{-N}=0$. Therefore, they are both nonzero (otherwise, they are both zero and again the dimension of the kernel will be at least 4).

Similarly, using the symmetry of the semiflow, we obtain that $d_{++}$ is a multiple of $d_{--}$, and $d_{+-}$ is a multiple of $d_{-+}$. Therefore, within the submatrix $\Delta_{0,1}$ all entries are dependent on one another, i.e., either they are all zero, or all nonzero. In case they are all zero, we obtain that the dimension of the kernel of $\Delta$ is 4, which is a contradiction. Therefore, they are all nonzero.

For the induction step, we have to argue how the connection matrix changes as the parameter $\lambda$ crosses a bifurcation value $\lambda_n$. Nevertheless, understanding changes of the connection matrix through one pitchfork bifurcation, as we did, is enough to replicate this in the proof of the general case. 

\subsection{Heteroclinic detection for general parabolic equations}

Consider
\begin{equation}\label{PDEquasi}
    u_t=u_{xx}+f(u)
\end{equation} 

The Conley index can be applied to detect heteroclinics as follows. Construct a closed neighborhood $N$ such that its maximal invariant subspace is the closure of the set of heteroclinics between $u_\pm$,
\begin{equation*}
    \Sigma=\{u_\pm\}\cup\overline{W^u(u_-)\cap W^s(u_+)}.
\end{equation*} 

This is also called liberalism in \cite{FiedlerRocha96}. Consider hyperbolic equilibria $u_-,u_+$ such that $i(u_-)=i(u_+)+1$ and satisfies both the Morse and the zero number permit conditions. Without loss of generality, assume $u_-(0)>u_+(0)$. 

It is used the Conley index to detect orbits between $u_-$ and $u_+$. Note that the semiflow generated by the equation \eqref{PDEquasi} on the Banach space X is admissible for the Conley index theory in the sense of \cite{Rybakowski82}, due to a compactness property that is satisfied by the parabolic equation \eqref{PDEquasi}, namely that trajectories are precompact in phase space. See Theorem 3.3.6 in \cite{Henry81}.

Suppose, towards a contradiction, that there are no heteroclinics connecting $u_-$ and $u_+$, that is, $\Sigma=\{ u_-,u_+\}$. Then, the index is given by the wedge sum $h(\Sigma)=[\mathbb{S}^n]\vee [\mathbb{S}^m]$, where $n,m$ are the respective Morse index of $u_-$ and $u_+$.

If, on the other hand, one can prove that $h(\Sigma)=[0]$, where $[0]$ means that the index is given by the homotopy equivalent class of a point. This would yield a contradiction and there should be a connection between $u_-$ and $u_+$. Moreover, the Morse-Smale structure excludes connection from $u_+$ to $u_-$, and hence there is a connection from $u_-$ to $u_+$.

Hence, there are three ingredients missing in the proof: the Conley index can be applied at all, the construction of a isolating neighborhood $N$ of $\Sigma$ and the proof that $h(\Sigma)=[0]$.

As mentioned above, in order to apply the Conley index concepts we need to construct appropriate neighborhoods and show that the Conley index is $[0]$. 

Consider the closed set
\begin{equation*}
    K({u_\pm}):= \left\{ u \in X \mid 
        \begin{array}{c} 
        z(u-u_-)=i(u_+)=z(u-u_+) \\
        u_+(0)\leq u(0)\leq u_-(0) 
        \end{array} \right\}
\end{equation*}

Consider also closed $\epsilon$-balls $B_\epsilon(u_\pm)$ centered at $u_\pm$ such that they do not have any other equilibria besides $u_\pm$, respectively, for some $\epsilon>0$. 

Define
\begin{equation*}
    N_\epsilon (u_\pm):=B_\epsilon(u_-)\cup B_\epsilon(u_+) \cup K({u_\pm}).    
\end{equation*}

The zero number blocking condition implies there are no equilibria in $K({u_\pm})$ besides possibly $u_-$ and $u_+$. Hence, $N_\epsilon (u_\pm)$ also has no equilibria besides $u_-$ and $u_+$. 

Denote $\Sigma$ the maximal invariant subset of $N_\epsilon$. We claim that $\Sigma$ is the set of the heteroclinics from $u_-$ to $u_+$ given by $\overline{W^u(u_-)\cap W^s(u_-)}$. 

On one hand, since $\Sigma$ is globally invariant, then it is contained in the attractor $\mathcal{A}$, which consists of equilibria and heteroclinics. Since there are no other equilibria in $N_\epsilon (u_\pm)$ besides $u_\pm$, then the only heteroclinics that can occur are between them.

On the other hand, the Theorem that conditions the number of zeros in stable/unstable manifolds implies that along a heteroclinic $u(t)\in\mathcal{H}$ the zero number satisfies $z^t(u-u_\pm)=i(u_+)$ for all time, since $i(u_-)=i(u_+)+1$. Therefore $u(t)\in K({u_\pm})$ and the closure of the orbit is contained in $N_\epsilon (u_\pm)$. Since the closure of the heteroclinic is invariant, it must be contained in $\Sigma$.

Lastly, it is proven that $h(\Sigma)=[0]$ in three steps, yielding the desired contradiction and the proof of the theorem. We modify the first and second step from \cite{FiedlerRocha96}, whereas the third remain the same. 

In the first step, a model is constructed displaying a saddle-node bifurcation with respect to a parameter $\mu$, for $n:=z(u_+-u_-)\in\mathbb{N}$ fixed,
\begin{equation}\label{prototype2}
    v_t=a(\xi,v,v_\xi)[v_{\xi\xi}+\lambda_nv]+ g_n(\mu,\xi,v,v_\xi)
\end{equation}
where $\xi\in [0,\pi]$ has Neumann boundary conditions, $\lambda_n=-n^2$ are the eigenvalues of the laplacian with $\cos(n\xi)$ as respective eigenfunctions, and 
\begin{equation*}
    g_n(\mu,\xi,v,v_\xi):=\left(v^2+\frac{1}{n^2}v^2_\xi-\mu\right)\cos(n\xi).
\end{equation*}

For $\mu>0$, a simple calculation shows that $v_\pm=\pm \sqrt{\mu}\cos(n\xi)$ are equilibria solutions of \eqref{prototype2} such that
\begin{equation}
    z(v_+-v_-)=n
\end{equation}
since the $n$ intersections of $v_-$ and $v_+$ will be at its $n$ zeroes. 

Moreover, those equilibria are hyperbolic for small $\mu>0$, such that $i(v_+)=n+1$ and $i(v_-)=n$. Indeed, parametrize the bifurcating branches by $\mu=s^2$ so that $v(s,\xi)=s\cos(n\xi)$, where $s>0$ correspond to $v_+$ and $s<0$ to $v_-$. Linearizing at the equilibrium $v(s,\xi)$ and noticing some terms cancel, the eigenvalue problem becomes
\begin{equation*}
    \eta v=a(\xi,s\cos(n\xi),sn\sin(n\xi))[v_{\xi\xi}+\lambda_n v]+\left[2s\cos(n\xi) v+\frac{2s\sin(n\xi) v_\xi}{-n}\right]\cos(n\xi)
\end{equation*}
where the unknown eigenfunction is $v$, corresponding to the eigenvalue $\eta$.

Hence $\eta_n(s)=2s$ is an eigenvalue with $v(s,\xi)$ its corresponding eigenfunction. 
Hence, by a perturbation argument in Sturm-Liouville theory, that is $\mu=0$ we have the usual laplacian with $n$ positive eigenvalues and one eigenvalue $\eta_n(0)=0$. Hence for small $\mu<0$, the number of positive eigenvalues persist, whereas for small $\mu>0$, the number of positive eigenvalues increases by 1. This yields the desired claim about hyperbolicity and the Morse index.

Now consider the quasilinear parabolic equation such that \eqref{prototype2} is the equilibria equation. The equilibria $v_\pm$ together with their connecting orbits form an isolated set
\begin{equation*}
    \Sigma_\mu(v_\pm):= \overline{W^u(v_-)\cap W^s(v_+)}
\end{equation*}
with isolating neighborhood $N_\epsilon(v_\pm)$, and the bifurcation parameter $\mu$ can also be seen as a homotopy parameter. Hence the Conley index is of a point by homotopy invariance as desired, that is,
\begin{equation}\label{conley1q}
    h(\Sigma_\mu(v_\pm))=h(\Sigma_0(v_\pm))=[0].
\end{equation}

In the second step, the $v_-$ and $v_+$ are transformed respectively into $u_-$ and $u_+$.

Recall $n=z(v_--v_+)=z(u_+-u_-)$. Hence, choose $\xi(x)$ a smooth diffeomorphism of $[0,\pi]$ that maps the zeros of $v_-(\xi)-v_+(\xi)$ to the zeros of $u_-(x)-u_+(x)$. Therefore, the zeros of $v_-(\xi(x))-v_+(\xi(x))$ and $u_-(x)-u_+(x)$ occur in the same points in the variable $x\in [0,\pi]$.

Consider the transformation 
\begin{align*}
    L: X &\to X   \\ 
    v(\xi)&\mapsto l(x)[v(\xi(x))-v_-(\xi(x))]+u_-(x)
\end{align*}
where $l(x)$ is defined pointwise through
\begin{align*}
    l(x):=
    \begin{cases}
        \frac{u_+(x)-u_-(x)}{v_+(\xi(x))-v_-(\xi(x))} &, \text{ if } v_+(\xi(x))\neq v_-(\xi(x)) \\
        \frac{\partial_x(u_{+}(x)-u_{-}(x))}{\partial_x(v_{+}(\xi(x))-v_{-}(\xi(x)))} &, \text{ if } v_+(\xi(x))= v_-(\xi(x))
    \end{cases}
\end{align*}
such that the coefficient $l(x)$ is smooth and nonzero due to the l'H\^opital rule. Hence, $L(v_-)=u_-$ and $L(v_+)=u_+$ as desired. Note we supposed $2\alpha+\beta>1$ so that solutions $u_\pm\in C^1$, hence $L$ is of this regularity as well. Moreover, $L$ is invertible with inverse having the same regularity. In particular, it is a homeomorphism, and hence a homotopy equivalence.

Moreover, the number of intersections of functions is invariant under the map $L$,
\begin{equation}
    z(L (v(\xi)-\tilde{v}(\xi)))=z(l(x)[v(\xi(x))-\tilde{v}(\xi(x))])=z(v(x)-\tilde{v}(x))
\end{equation}
and hence $K({v_\pm})$ is mapped to $K({u_\pm})$ under $L$. 

Consider $w(t,x):=L(v(t,\xi))$, hence the map $L$ modifies the equation \eqref{prototype2} into the following equation
\begin{equation}\label{IDK}
    w_t=\tilde{a}(x,w,w_x)w_{xx}+\tilde{f}(x,w,w_x)
\end{equation}
where the Neumann boundary conditions are preserved, and the terms $\tilde{a},\tilde{f}$ are
\begin{align*}
    \tilde{a}(x,w,w_x):=&\frac{x_\xi^2}{l(x)} \cdot a(x,L^{-1}(w),\partial_xL^{-1}(w))\\
    \tilde{f}(x,w,w_x):=&g_n(\mu,x,L^{-1}(w),\partial_xL^{-1}(w))+ (w_x\cdot x_{\xi\xi}-\partial_\xi^2u_-)-\frac{l_{\xi\xi} \cdot (w-u_-)_\xi}{l}\\
    &-\lambda a(x,L^{-1}(w),\partial_xL^{-1}(w))\cdot L^{-1}(w).
\end{align*}

Note that the equilibria $v_\pm$ are mapped into $w_\pm:=L(v_\pm)=u_\pm$, which are equilibria of \eqref{IDK}, with same zero numbers and Morse indices as $v_\pm$ and $u_\pm$.

The isolated invariant set $\Sigma_\mu(v_\pm)$ is transformed into $L(\Sigma_\mu(v_\pm))=\Sigma_\mu(w_\pm)$, which is still isolated and invariant, with invariant neighborhood $L(N_\epsilon(v_\pm))=N_\epsilon(w_\pm)$. Moreover, the Conley index is preserved, since $L$ is a homotopy equivalence,
\begin{equation}\label{conley2q}
   h(\Sigma_\mu(v_\pm))= h(L(\Sigma_\mu(v_\pm)))=h(\Sigma_\mu(w_\pm)).
\end{equation}

Hence, one identifies the equilibria $v_\pm$ in the model constructed \eqref{prototype2} with the equilibria $w_\pm=u_\pm$ from the equation \eqref{PDEquasi}, by preserving neighborhoods and the Conley index, since $L$ is a homotopy equivalence. The identified equilibria $u_\pm$ satisfy the equation \eqref{IDK}, and we still have to modify it to become \eqref{PDEquasi}. For such, we perform now a last homotopy between the solutions $w$ and $u$.

In the third step, we homotope the diffusion coefficient $\tilde{a}$ and nonlinearity $\tilde{f}$ from the equation \eqref{IDK} to be the desired diffusion $a$ and reaction $f$ from the equation \eqref{PDEquasi}. Indeed, consider the parabolic equation 
\begin{equation*}
    u_t=a^\tau(x,u,u_x)u_{xx}+f^\tau(x,u,u_x)
\end{equation*}
where
\begin{align*}
    a^\tau&:=\tau \tilde{a}+(1-\tau)a+\sum_{i=- \text{ , } +}\chi_{u_i}\mu_{u_i}(\tau)[u-{u_i}(x)]\\
    f^\tau(x,u,u_x)&:=\tau \tilde{f}+(1-\tau)f+\sum_{i=- \text{ , } +}\chi_{u_i}\mu_{u_i}(\tau)[u-{u_i}(x)]
\end{align*}
and $\chi_{u_i}$ are cut-offs being 1 nearby $u_i$ and zero far away, the coefficients $\mu_i(\tau)$ are zero near $\tau=0$ and $1$ and shift the spectra of the linearization at $u_\pm$ such that uniform hyperbolicity of these equilibria is guaranteed during the homotopy. Note that $u_\pm$ have the same Morse indices, as solutions of both equations \eqref{PDEquasi} and \eqref{IDK}. Therefore, the $\mu_i(\tau)$ only makes sure none of these eigenvalues cross the imaginary axis. 

Consider $u_\pm$ and their connecting orbits during this homotopy, 
\begin{equation*}
    \Sigma^\tau(u_\pm):=\overline{W^u(u_-)\cap W^u(u_+)}.
\end{equation*}

Note that $\Sigma^\tau(u_\pm)\subseteq K({u_\pm})$, for all $\tau\in [0,1]$, since the dropping lemma holds throughout the homotopy. The equilibria $u_\pm$ do not bifurcate as $\tau$ changes, due to hyperbolicity. Choosing $\epsilon>0$ small enough, the neighborhoods $N_\epsilon(u_\pm)$ form an isolating neighborhood of $\Sigma^\tau(u_\pm)$ throughout the homotopy. Indeed, $\Sigma^\tau(u_\pm)$ can never touch the boundary of $K({u_\pm})$, except at the points $u_\pm$ by the dropping lemma. Once again the Conley index is preserved by homotopy invariance,
\begin{equation}\label{conley3q}
    h(\Sigma(u_\pm))=h(\Sigma^0(u_\pm))=h(\Sigma^\tau(u_\pm))=h(\Sigma^1(u_\pm))=h(\Sigma_\mu(w_\pm)).
\end{equation}

Finally, the equations \eqref{conley1q}, \eqref{conley2q} and \eqref{conley3q} yield that the Conley index of $\Sigma$ is the homotopy type of a point, and hence the desired result:
\begin{equation}\label{conley4q}
    h(\Sigma(u_\pm))=h(\Sigma_\mu(w_\pm))=h(\Sigma_\mu(v_\pm))=[0].
\end{equation}

\end{document}